\documentclass[12pt, a4paper,twoside]{paper}
\usepackage[utf8]{inputenc}
\usepackage[T1]{fontenc}
\usepackage[english]{babel}
\usepackage{lmodern}
\usepackage{newlfont}
\usepackage{geometry}
\usepackage{epigraph}
\usepackage[utf8]{inputenc}
\usepackage[T1]{fontenc}
\usepackage[english]{babel}
\usepackage{lmodern}
\usepackage{newlfont}
\usepackage{epigraph}
\usepackage{hyperref}

\usepackage{epigraph}
\usepackage{csquotes}
\usepackage{tikz-cd}
\usepackage{graphicx}

\usepackage{titlesec}
\usepackage{color}
\usepackage{blindtext}
\usepackage{wrapfig}
\usetikzlibrary{arrows.meta, patterns}

\usepackage{bm}
\usepackage[font=scriptsize,format=plain, labelfont=bf, ]{caption}
\usepackage{tabularx}
\usepackage{stackengine}
\usepackage{dsfont}
\usepackage[normalem]{ulem}
\usepackage{mathtools}
\usepackage{amsmath}
\usepackage{amssymb}
\usepackage{mathrsfs}
\usepackage{amsthm}
\usepackage{accents}

\usepackage[sorting=none]{biblatex}

\usepackage{wrapfig}

\usepackage{csquotes}
\usepackage{tikz-cd}
\usepackage{graphicx}
\usepackage{titlesec}
\usepackage{color}
\usepackage{blindtext}
\usepackage{wrapfig}
\usepackage{bm}
\usepackage{tabularx}
\usepackage{stackengine}
\usepackage{dsfont}

\usepackage{mathtools}
\usepackage{amsmath}
\usepackage{amssymb}
\usepackage{mathrsfs}
\usepackage{amsthm}
\usepackage{accents}

\usepackage[sorting=none]{biblatex}
\usepackage[title]{appendix}

\theoremstyle{definition}

\newtheorem{proposition}{Proposition}[section]

\newtheorem{definition}[proposition]{Definition}

\newcommand{\beq}{\begin{equation}}
\newcommand{\eeq}{\end{equation}}

\usepackage{lineno}
\begin{document}
 
%
\title{Data-Driven Reduced Modeling of Recurrent Neural Networks}

\maketitle
\centerline{A. Marraffa\footnote{Institute for Mechanical Systems, ETH Zurich, Zurich, Switzerland,\label{note1}}, R. Krause\footnote{Institute of Neuroinformatics, University of Zurich \& ETH Zurich, Zurich, Switzerland\label{note2}}\footnote{Neuroscience Center Zurich, University of Zurich \& ETH Zurich, Zurich, Switzerland\label{note3}}, V. Mante\footref{note2}\footref{note3} and G. Haller\footref{note1}}






\medskip
\centerline{\large \bf  
}

\medskip
\centerline{\large \bf } 

\bigskip

\begin{abstract}




Artificial Recurrent Neural Networks (RNNs) are widely used in neuroscience to model the collective activity of neurons during behavioral tasks. The high dimensionality of their parameter and activity spaces, however, often make it challenging to infer and interpret the fundamental features of their dynamics.

In this study, we employ recent nonlinear dynamical system techniques to uncover the core dynamics of several RNNs used in contemporary neuroscience. Specifically, using a data-driven approach, we identify Spectral Submanifolds (SSMs), i.e., low-dimensional attracting invariant manifolds tangent to the eigenspaces of fixed points. The internal dynamics of SSMs serve as nonlinear models that reduce the dimensionality of the full RNNs by orders of magnitude. 

Through low-dimensional, SSM-reduced models, we give mathematically precise definitions of line and ring attractors, which are intuitive concepts commonly used to explain decision-making and working memory. The new level of understanding of RNNs obtained from SSM reduction enables the interpretation of mathematically well-defined and robust structures in neuronal dynamics, leading to novel predictions about the neural computations underlying behavior.

\end{abstract}


    \section{Introduction}

Behavior emerges from concerted activity across large populations of neurons. In the study of neural population activity, two primary questions emerge: How do neurons collectively create an internal representation of the external world based on the inputs they receive, and how do neural computations on these representations lead to behavior?

Recent advances in neuronal modeling and analysis have made it possible to address these questions, providing new insights into the nature of computations implemented at the level of neural populations. To establish a connection between the firing patterns of neural populations and low-dimensional, interpretable behavioral outputs, individual neurons can be interpreted as degrees of freedom within a dynamical system. Computations can then be understood in terms of the evolution of this dynamical system when driven by the inputs it receives from the senses or from other areas \cite{computation_through_dynamics}, \cite{BrodyOpinion}. Dynamical systems modeling naturally lends itself to analyses based on systems identification, whereby neural computations reflect the properties and interaction of low-dimensional, attracting structures, even when single-neuron responses are not directly interpretable.

Accumulating evidence from such approaches suggests that specific classes of dynamical motifs are consistently involved in cognitive and motor operations. For instance, saddle-point dynamics have been associated with decision-making processes \cite{MachensBrody}, \cite{WangDecisionMI}, and rotational dynamics have been proposed as a hallmark of population activity underlying movement preparation and execution \cite{Shenoy2013Movement} \cite{OscillReaching}. Some other heuristics are also helpful in characterizing observed dynamics in more complex tasks: "ring attractors" have been proposed as substrates for working memory \cite{nhimsrnn}, \cite{ostojic2023} \cite{memorychad} \cite{Khona2021AttractorAI}, \cite{wm2}, \cite{wm3}, and "line attractors" have been discussed for tasks that require the accumulation of evidence \cite{mante_context-dependent_2013}\cite{la1}, \cite{la2}, \cite{wm1}. Importantly, more involved and structured behavioral tasks could be understood as the composition of such dynamical primitives \cite{Pouget}, \cite{ringattractor}.

Identifying such dynamical motifs requires accurate models of neural dynamics. To this end, a wide range of methods has been developed in recent years to infer dynamics from population—level recordings. Many techniques rely on fitting a model of dynamics to the data (see \cite{GPFA} for an exception), but generally suffer from a tradeoff between complexity, expressive power, and the interpretability of the fitted models.
At one extreme, machine—learning—based methods such as LFADS \cite{lfadssussillo} and CEBRA \cite{CEBRA} (see also \cite{BrodyRecent}, \cite{MARBLE}) can be used to fit highly nonlinear, even chaotic dynamics, but often at the expense of mechanistic insight and predictive power. At the other end are models based on linear dynamical systems, either a single one \cite{sahani_early}, or switching \cite{lind_switching}, piecewise linear \cite{LindermanSussillo}, or time-dependent linear systems \cite{Galgali2021ResidualDR}. Since the possible dynamics of linear systems are well understood, the resulting descriptions of dynamics are highly interpretable, but cannot capture intrinsically nonlinear behavior, such as coexisting attractors and chaotic dynamics.

Here, we propose a new approach to estimate neural dynamics based on a recent nonlinear model reduction technique that strikes a balance between the ability to capture intrinsically nonlinear dynamics and interpretability. Central to our approach is the recent theory of Spectral Submanifolds (SSMs)\cite{firstSSM}, which are ubiquitous and robust low-dimensional, attracting invariant manifolds that emanate from steady states with nondegenerate and nonresonant linearized spectra (see \cite{SSMReview} for a general introduction, applications, and further references). The internal dynamics of SSMs serve as low-dimensional, explicit, polynomial models that explain the dominant slow nonlinear dynamics of high-dimensional dynamical systems.
Based on these reduced models, dynamical motifs in neural networks can be analyzed and explained in a mathematically rigorous way.

We demonstrate the power of data—driven SSM—based modeling for neural dynamics using simulation output data from recurrent neural networks (RNNs). RNNs trained on specific behavioral tasks have emerged as a key tool to reproduce the features of neural dynamics in brain recordings and to generate hypotheses about the neural computations implemented by the underlying neural circuits \cite{RNNRev}, \cite{RobustRNN}, \cite{Gerstner}. However, obtaining a complete understanding of the dynamics of RNNs is often challenging, even when all their underlying parameters are fully known. Commonly used techniques for identifying dynamical motifs in RNNs rely on the search for fixed points and "slow points" \cite{slowpoints}, \cite{mante_context-dependent_2013}, \cite{ringattractor}, i.e., locations in neural state space where dynamics does not change or varies slowly. The geometry of slow points is largely conserved across different RNNs solving the same tasks and is therefore considered a useful heuristic of the computations they implement \cite{sussUniversality}. Further insights into computations can be obtained through linear approximations of the dynamics around the identified fixed points. 

Such standard approaches to understanding RNN computations have several downsides, which are overcome by SSM-based modeling. First, slow points are typically found numerically based on subjectively chosen thresholds on the rate of change of neural dynamics. Second, linear approximations of non-linear systems are valid only only locally around fixed points, but not around slow points. Third, even when dynamics can be linearized, around hyperbolic fixed points the linear approximations are only valid in domains with simple, \textit{linearizable} dynamics. These domains cannot contain multiple isolated fixed points, periodic orbits, quasiperiodic tori, chaotic attractors, or transitions among such invariant sets. In contrast, we show that data-driven algorithms to estimate low-dimensional SSM and their reduced dynamics (\cite{DataDriven}) can provide a complete, yet simple description of both linear and nonlinear dynamics across all computationally relevant activity domains. The resulting mathematical descriptions of the global dynamics provide the basis for more precise hypotheses about the nature of neural computations in biological neural networks.

    \section{Results}

We model recurrent neural networks with evolution rules of the form: \begin{equation}\label{eq:RNNeq}
 \frac{dx}{dt} = R_{\theta} (x(t), u(t)) 
\end{equation} where $x(t) \in \mathbb{R}^N$ is the network state variable, $u(t) \in\mathbb{R}^{N_{in}} $ is the input vector and $R_{\theta}$ is a function depending on the parameter vector $\theta$ that defines the network dynamics.

Our approach to analyzing autonomous, high-dimensional RNN dynamics begins by assuming the existence of a known steady-state solution, which in practice can be identified as zeros of the right-hand-side of Eq. \ref{eq:RNNeq}. Around such a steady state, we can deduce the existence of attracting invariant SSMs.

SSMs can be technically defined as the smoothest invariant manifolds that are tangent to dominant eigenspaces of steady states. The existence of these manifolds is ensured under nonresonance conditions on the spectrum of the linearized dynamical system at the steady state \cite{SSMReview}. Once the existence of an SSM is established, it can be approximated through a polynomial expansion over the dominant spectral subspace used in its identification. While the existence of SSMs is inferred from the linearized spectrum at a steady state, these invariant manifolds do extend beyond the domain of linearization and hence can carry characteristically nonlinear dynamics. The internal dynamics of attracting SSMs serve as a mathematically exact reduced-order model with which nearby trajectories synchronize exponentially fast. 

We employ the SSMLearn algorithm \cite{DataDriven} to find explicit polynomial parametrizations for low-dimensional, attracting, invariant SSMs attached to steady-state solutions. SSMLearn then approximates the reduced polynomial ordinary differential equation (ODE) that governs the system evolution restricted to the SSM, providing a low-dimensional, polynomial model that captures the core dynamics of the original RNN.

A key advantage of this approach is the substantial reduction in model dimensionality, which converts a high-dimensional RNN into a low-dimensional polynomial ODE. Therefore, SSM-reduced models allow for the explicit identification of phase space structures, such as fixed points, limit cycles, and bifurcations, that determine the global asymptotic behavior of the network. To assess the effect of system parameters on dynamics \ref{eq:RNNeq}, we can construct a slow manifold for each value of a given system parameter and view it as a section of a global slow manifold in the extended phase space that includes the system parameters as well. This, in turn, enables us to build parameter-dependent polynomial models for the SSM-reduced dynamics.

\subsection{A Model for Context-Dependent Decision-Making}
  \begin{figure}
\centering
\includegraphics[width = .48\textwidth]{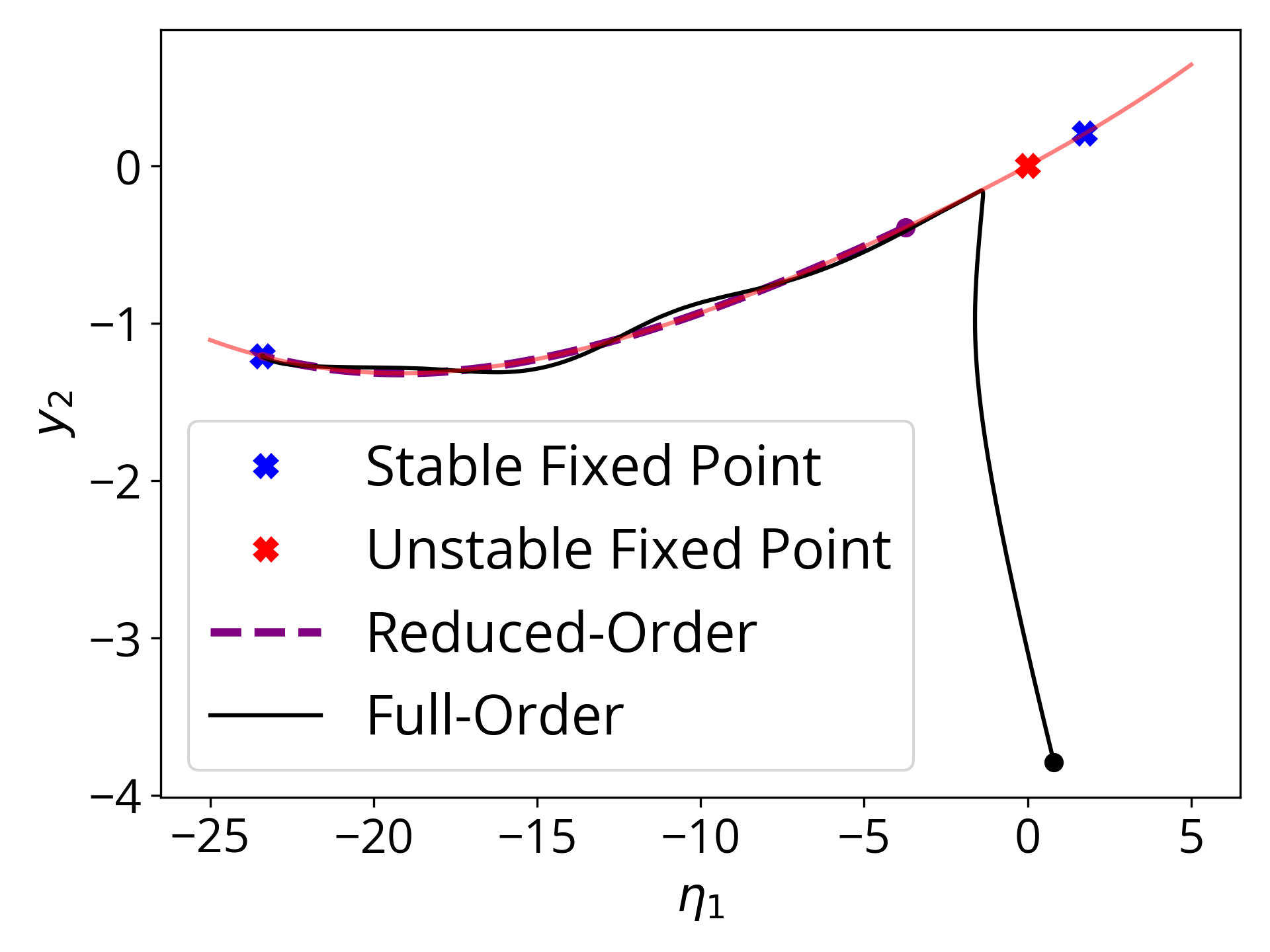}
\includegraphics[width = .48\textwidth]{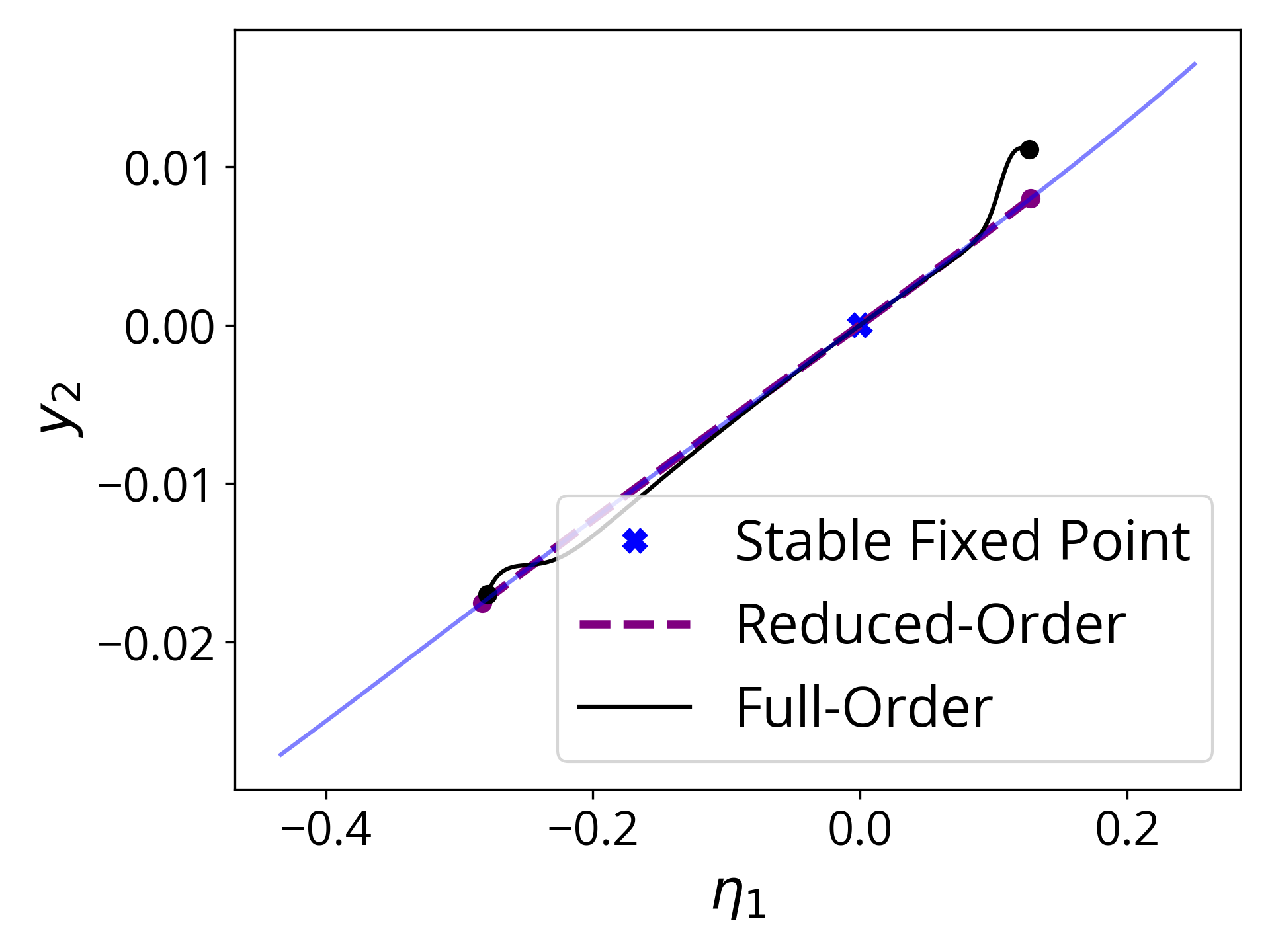}
\centering
\includegraphics[width = .48\linewidth]{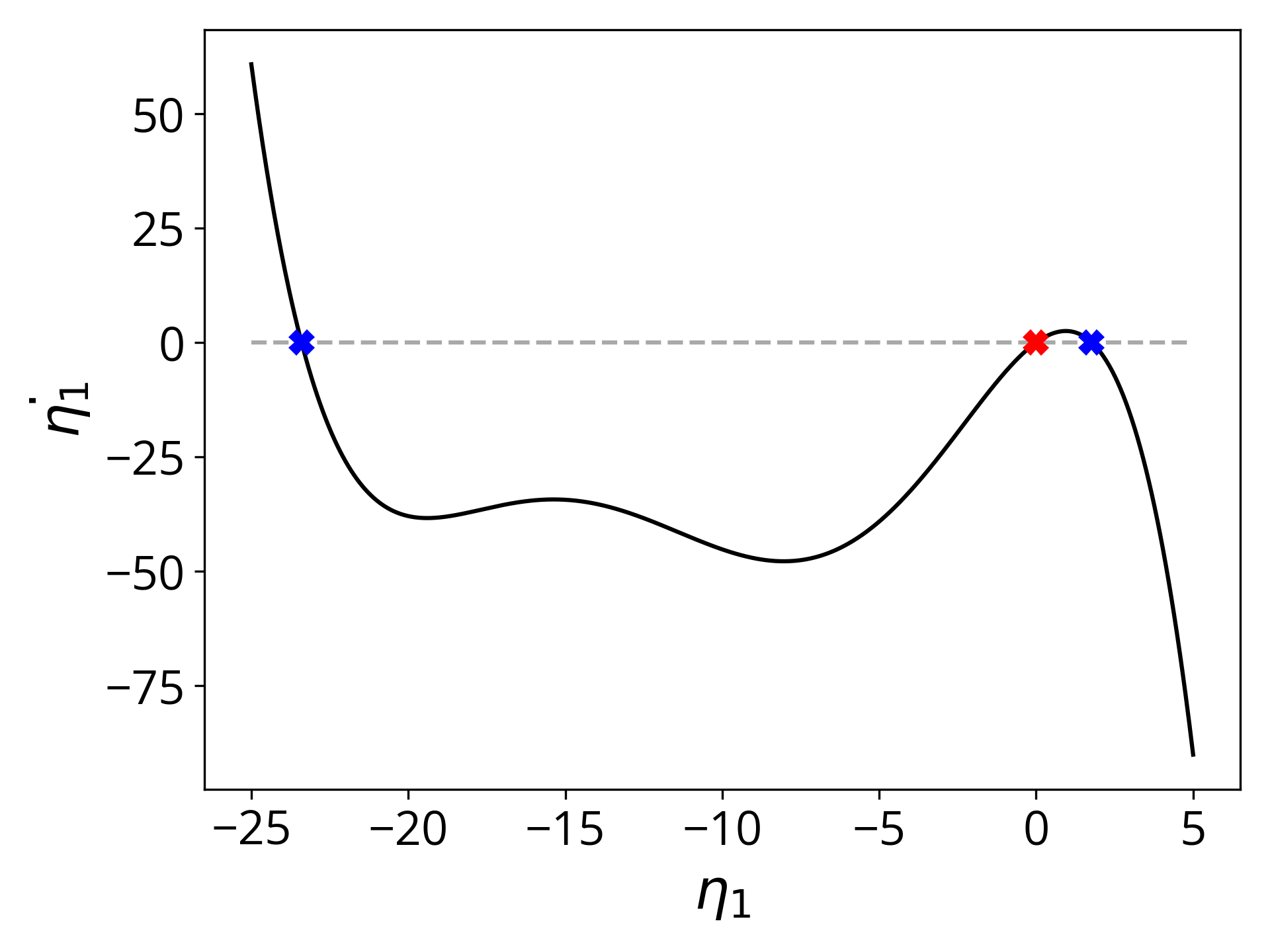}
\includegraphics[width = .48\linewidth]{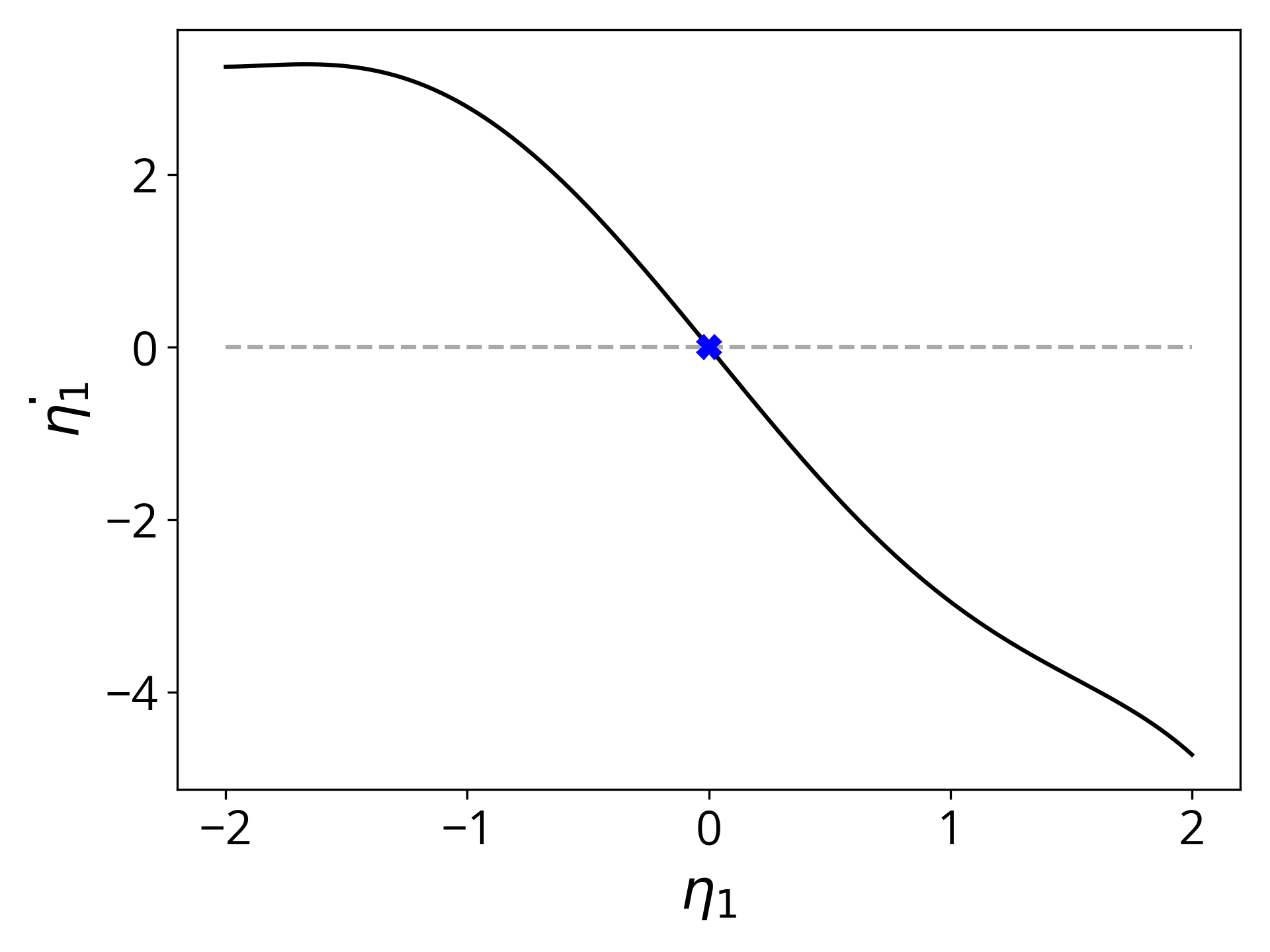}
\centering
\includegraphics[width = .48\linewidth]{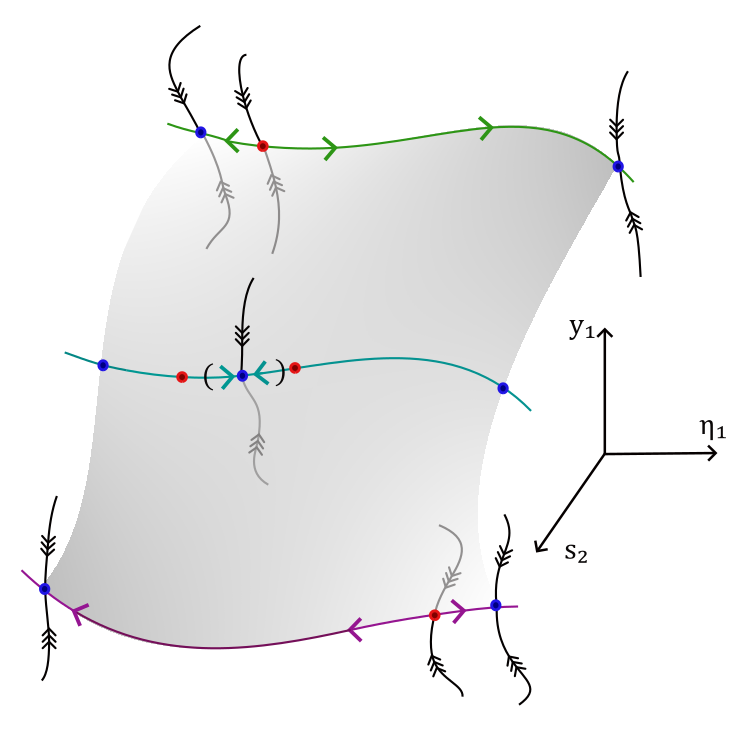}
\centering
\includegraphics[width = .48\linewidth]{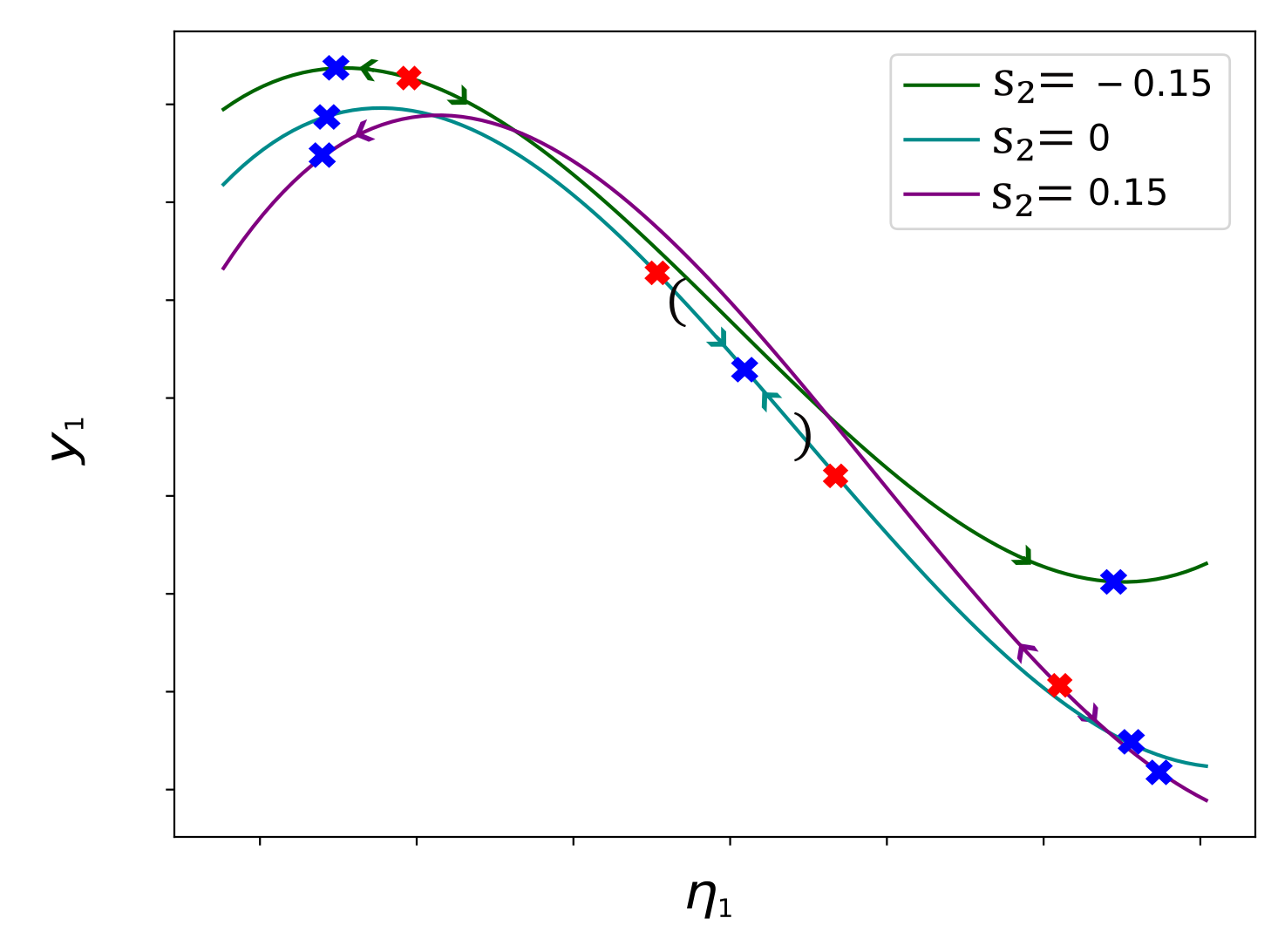}
\caption{SSMs carrying the RNN reduced dynamics. (Upper left) Unstable manifold at order five in coordinates $(\eta_1, y_2)$, where $\eta_1$ parametrizes the spectral subspace $E_1$ and $y_i = x_i-x_{i,0}$ are the original RNN coordinates centered around the unstable fixed point. We depict a full-order test trajectory and the corresponding reduced trajectory on the SSM (after $t\approx\frac 1  {\lambda_2}$, when transients have decayed, purple dot.) converging to the fixed point with the larger domain of attraction along the unstable manifold. Black dots correspond to initial conditions. The fixed points in the plot (crosses) are the fixed points of the full model, lying on the unstable manifold and correctly reproduced by the reduced-order model. The Manifold Fitting Error (MFE), the mean distance between observed trajectories and trajectories projected on the manifold, has magnitude $\approx 0.007$, and the Normal Mean Trajectory Error (NMTE), the mean reconstruction error of the reduced model, is $\approx 0.03$. (Upper right) The slowest SSM at order three and test full-order and reduced-order trajectories plotted in coordinates $(\eta_1, y_2)$, with $y_i$ now centered around the stable fixed point. The MTE and NMTE have magnitude $\approx 0.05$. (Center) Graphs of the right-hand sides of the reduced-order models: on the 1D unstable manifold (left) and on the 1D slowest SSM (right), both at order five. On the left plot, the unstable fixed point (red cross) serves as a boundary for the domains of attraction of the two stable fixed points (blue crosses) on the unstable manifold. On the right plot, we observe convergence to the stable fixed point. (Bottom) Sketch of the robust (normally hyperbolic) 2D invariant manifold (\textit{slow manifold)} in the extended phase space constructed including the $s_2$-parameter direction (left), in coordinates $(\eta_1, s_2, y_1)$, where $\eta_1$ and $y_1$ are as above, and in coordinates $(\eta_1, y_1)$ for selected $s_2$ values (right). We see how the 1D latent dynamics underlying flexible decision-making change when varying the relevant sensory input $s_2$ from negative to positive values, passing through 0, resulting in a change in the preferred choice for the network. When $s_2=0$, there is a stable fixed point that attracts initial conditions around zero and, locally, the slow manifold coincides with its 1D slowest SSM (in round brackets).}
\label{fig:unstablemfld}
\end{figure}


As a first example, we analyze the $100$D RNN performing the context-dependent decision-making task introduced by \cite{mante_context-dependent_2013}, as described in the Methods section \ref{ss:cddmRNN}. This network receives as input two sensory cues, $s_1$ and $s_2$, whose values can be positive or negative, and two context cues, $(c_1, c_2)$ with values in $\{(0,1), (1,0)\}$. The network's output is expressed through a scalar readout function that will have the same sign of $s_1$, if $(c_1,c_2)=(1,0)$, or $s_2$, $(c_1,c_2)=(0,1)$.

The authors of \cite{mante_context-dependent_2013} explain the network dynamics underlying flexible decision-making using the concept of a \textit{line attractor}, which they view as a curve in phase space consisting of \textit{slow points} \cite{badLinearization}, \cite{mante_context-dependent_2013}, \cite{slowpoints}. Although not an attractor in a strict mathematical sense, the line attractor is a heuristic used to explain what is observed in the first principal component space: trajectories are attracted to this curve when sensory inputs are off and diverge from it once sensory inputs are turned on.

Within this 1D structure, \cite{mante_context-dependent_2013} finds two stable fixed points that represent the two possible decisions in the task. The final state of the network, i.e., the fixed point to which it converges, is then argued to depend on the direction of perturbation away from the line attractor, determined by the direction and magnitude of the input vector (see Fig. 3 in \cite{mante_context-dependent_2013} and Fig. 2 in \cite{pagan}).

Here we revisit this interpretation through SSM theory, reinterpreting the RNN behavior as reduced dynamics on a 1D invariant slow manifold. 
We begin by analyzing the system under fixed context and sensory input, considering two conditions: the sensory input-off regime, corresponding to the cue period before sensory integration in animal experiments, and the sensory input-on regime, where active sensory processing occurs.

In both cases, we identify at least one fixed point in the system's phase space, which ensures the existence of SSMs for the flow map.
In the sensory input-on case, the network has an unstable fixed point and two stable fixed points, each corresponding to a distinct choice in the task. Importantly, these two stable fixed points are present for every input parameter value in the range explored in this experiment (see Fig. \ref{fig:bif_diag} in the Supplementary Materials). Linearization around the unstable fixed point reveals a single unstable eigenvalue, indicating a 1D unstable subspace. We use SSMLearn to find the expansion coefficients of the corresponding unstable manifold of the nonlinear system (Fig. \ref{fig:unstablemfld}, upper left and Fig. \ref{fig:SupSSMs}, right) up to order five and derive the associated reduced-order model
\begin{equation}
    \dot{\eta}_1 = \sum_{i=1}^{5} a_i \, \eta_1^i + \mathcal{O}(|\eta_1|^6)
\label{eq:rom1Dunstmfld}
\end{equation}
 whose right-hand side is plotted in Fig. \ref{fig:unstablemfld} (center left). The SSM-reduced model (\ref{eq:rom1Dunstmfld}) at this order of approximation enables accurate prediction of the full RNN dynamics from simulations of the reduced model.

In fact, the 1D model (\ref{eq:rom1Dunstmfld}) is sufficient to reproduce the RNN behavior during decision making. On the unstable manifold, the size of the domains of attraction of the two stable fixed points determines the global phase space geometry: random initial conditions are more likely to evolve toward the fixed point associated with the larger domain of attraction, which corresponds to the correct choice in the task. We estimate these domains of attraction in the full phase space using Finite-Time Lyapunov Exponent (FTLE) analysis in the Supplementary Materials \ref{appC}. FTLE ridges identify regions of phase space where initial conditions close to each other have a high rate of separation. In this context, these initial conditions mark the boundaries of the domains of attraction of the two fixed points (see Fig. \ref{fig:ftle2d}).


 In the sensory input-off case, the dynamics in a small neighborhood of the origin (where initial conditions are set in the experiments) are governed by a stable fixed point with one slowest eigenmode in its linearized dynamics. Using SSMLearn, we find the expansion coefficients the corresponding 1D slowest SSM (see Fig. \ref{fig:unstablemfld}, upper right). Full details and the corresponding equations for the reduced dynamics, whose graph is shown in Fig. \ref{fig:unstablemfld} (center right), are provided in the Supplementary Materials \ref{appD:const}.

These results lead to a re-interpretation of the line attractor in \cite{mante_context-dependent_2013}. The set of slow points identified in \cite{mante_context-dependent_2013} during the sensory input off-period corresponds to the slowest SSM anchored to the stable fixed point. When sensory inputs are off, this manifold attracts nearby trajectories and carries stable, near-equilibrium dynamics. These dynamics are perceived as slow once the trajectories are close enough to the stable fixed point, to which they converge in infinite time.

When sensory inputs are activated, the dynamics transition: trajectories converge to the decision-related stable fixed point along the 1D unstable manifold of a now unstable fixed point (in Fig. 3 of \cite{mante_context-dependent_2013}, trajectories labeled as sensory input on). The reduced dynamics restricted to this SSM exhibit a bistable structure, wherein the unstable fixed point acts as a boundary separating the domains of attraction of the two stable fixed points.

Importantly, as depicted in Fig. \ref{fig:unstablemfld} (bottom), a 1D slow manifold is present throughout the explored range of sensory inputs and carries the relevant dynamics throughout the decision-making process. Specifically, when sensory inputs are zero, the slow manifold coincides locally with the slowest SSM attached to a stable fixed point. As the relevant sensory input \( s_i \) (with \( i = 1 \) for context 1 and \( i = 2 \) for context 2) is varied from zero, we retrieve the slow manifold as the unstable manifold of a (different for different $s_i$ signs) unstable fixed point that separates the basins of attraction of the two stable equilibria (see Fig. \ref{fig:bif_diag} and discussion therein in the Supplementary Materials).

As we show in the Supplementary Methods \ref{appA:paramSSM}, we can assess the robustness of the 1D slow manifold as sensory inputs change, and we know that its geometry and associated reduced dynamics deform smoothly. In other words, in the extended phase space (constructed including the parameter axis, as in Fig. \ref{fig:unstablemfld}, bottom left), we find a 2D, at least $C^1$ attracting slow manifold that carries the latent dynamics underlying flexible decision-making. Around persisting fixed points, we can explicitly determine the dependencies on parameters through parameter-dependent polynomial expansions on their SSMs (see Supplementary Materials \ref{appD:par} for a complete discussion).

\begin{figure}
    \centering
 \includegraphics[width = \linewidth]{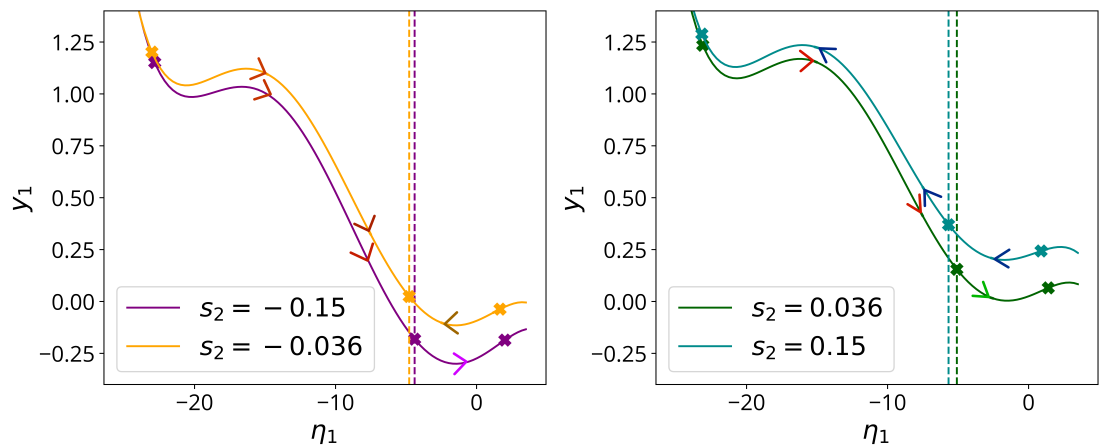}   
    \caption{One-dimensional slow manifolds (unstable manifolds) of the system for different values of the sensory input parameter $s_2$ (selected by context 2). The fixed point (coloured crosses) the network has learned to converge to depends on the sign of $s_2$. The way the RNN realizes this input-output relationship is changing location to the unstable fixed point (central cross) along the slow manifold depending on $s_2$: the unstable fixed point separates the domains of attraction of the two stable fixed points in the reduced dynamics, so its location (dashed line) determines the asymptotic behavior of trajectories, given their initial conditions. The crosses correspond to the fixed points, the two stable one divided by the unstable one. We display the arrows indicating effect of the input vector on the reduced dynamics (projection onto the SSM). This picture shows that the projection of the input vector (or its mean, in the noisy case) on the reduced dynamics does not determine the asymptotic behavior (choice) of the RNN, as we can see that, for example, it would lead to the wrong direction for $s_2= 0.036$ (red arrows). However, whenever the perturbed initial conditions lie within the correct domain of attraction, the convergence to the correct choice is guaranteed (see also \ref{appC}). }
    \label{fig:inputvec}
\end{figure}

Our analysis also clarifies the role of the input vector and the recurrent dynamics in determining the outcome of the decision-making task, allowing us in particular to distinguish among the three hypotheses proposed by \cite{pagan}. Notably, we find that the direction of the input vector in the reduced dynamics does not, by itself, determine the result of sensory integration (see Fig. \ref{fig:inputvec}). Moreover, as illustrated in Figs. \ref{fig:unstablemfld} (bottom), \ref{fig:inputvec}, and \ref{fig:kmflds}, the geometry of the SMMs remains largely constant across context and sensory inputs. As the input vector changes orientation, the tangent spaces of the SSMs do not align with it to facilitate integration.  Instead, convergence is determined by the location of initial conditions relative to the two domains of attraction. The width of these domains (and the location of the separating unstable fixed point on the slow manifold) varies parametrically with the relevant sensory input (or with its mean in the noisy case) to achieve the desired outcome. In other words, mean input vectors are parameters that determine the geometry of the phase space of the RNN.

In summary, based on the above findings we propose an alternative to the \textit{line attractor interpretation} of these RNN, namely a mathematically rigorous 1D \textit{slow manifold explanation}, and provide a method to systematically find these slow manifolds emanating from known steady-state solutions.
 
More generally, 1D slow manifolds may play a central role in decision-making tasks, not only in RNNs but also in biological neural systems. For instance, \cite{chadwick_mice} shows that a linear Vanilla RNN trained on mouse behavioral data from a sensory discrimination task learns a 1D slowest mode during training. This observation suggests that SSM theory has the potential to aid the analysis of neural computations in both artificial and biological settings. Importantly, the derivation of SSMs extends to cases where activity is corrupted by noise. We discuss the effects of adding noise to the RNN dynamics in the Supplementary Materials \ref{appD:const}, where we treat uniformly bounded Gaussian noise as a small-amplitude, time-dependent perturbation to the RNN studied so far.

\subsection{A Model for Oscillations with Input-Dependent Frequency}
Our second example, the $100$D \textit{Sine Wave Generator} RNN in \cite{renate} (see Method section \ref{ss:swgRNN}), is trained to produce a sinusoidal scalar readout $z(t)$, with oscillation frequency selected by a scalar input $u$ (see Eq. (\ref{eq:VanRNN})). 

In the phase space of the constant-input RNN, there is an unstable fixed point with two complex-conjugate unstable eigenvalues widely separated in their real parts from the stable ones. As a consequence, a 2D slow SSM exists in this problem and coincides with the unstable manifold of the fixed point. This SSM carries the observable oscillatory dynamics of the system. We find the SSM equations and its reduced dynamics with SSMLearn up to order three and six, which are the lowest expansion orders that give low invariance error (distance of trajectories from the SSM) and prediction error (distance between full and reduced trajectories), respectively. Specifically, the SSM-reduced dynamics takes the form
\begin{align}
\dot{\eta}_\ell &= \sum_{i=0}^{6} \sum_{j=0}^{6-i} a_{ij}^{(\ell)} \, \eta_1^i \eta_2^j \ + \ \mathcal{O}(7), \quad \ell = 1,2.
\label{eq:rom2Dswg}
\end{align}

Analyzing the SSM-reduced dynamics (\ref{eq:rom2Dswg}), we find a globally attracting limit cycle that is responsible for the periodic oscillations with the desired frequency. This type of dynamics is exactly the simplest possible realization of a dynamical system that oscillates with a single frequency. 

\begin{figure}
    \centering
    \includegraphics[width=0.48\linewidth]{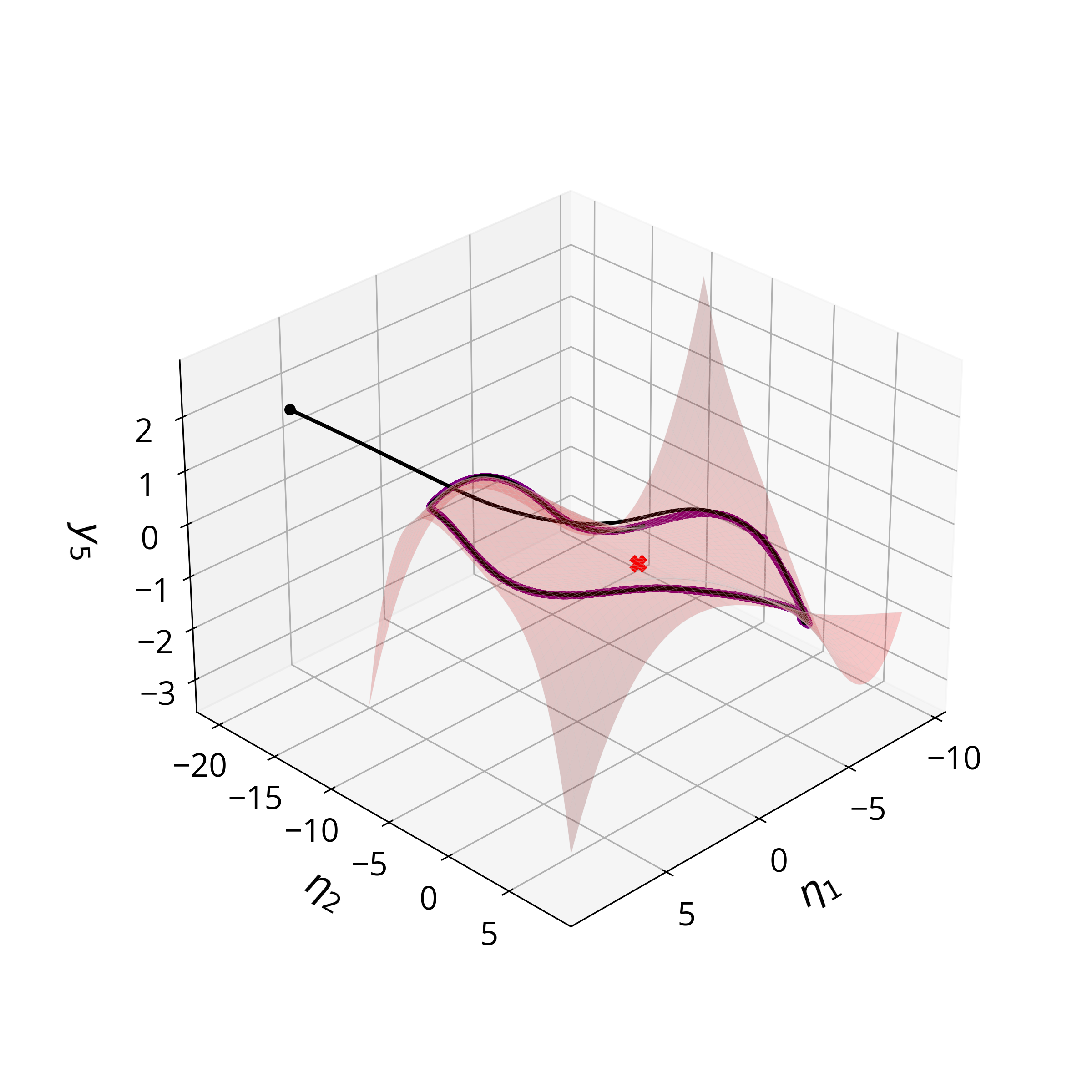}
    \centering
    \includegraphics[width=0.48\linewidth]{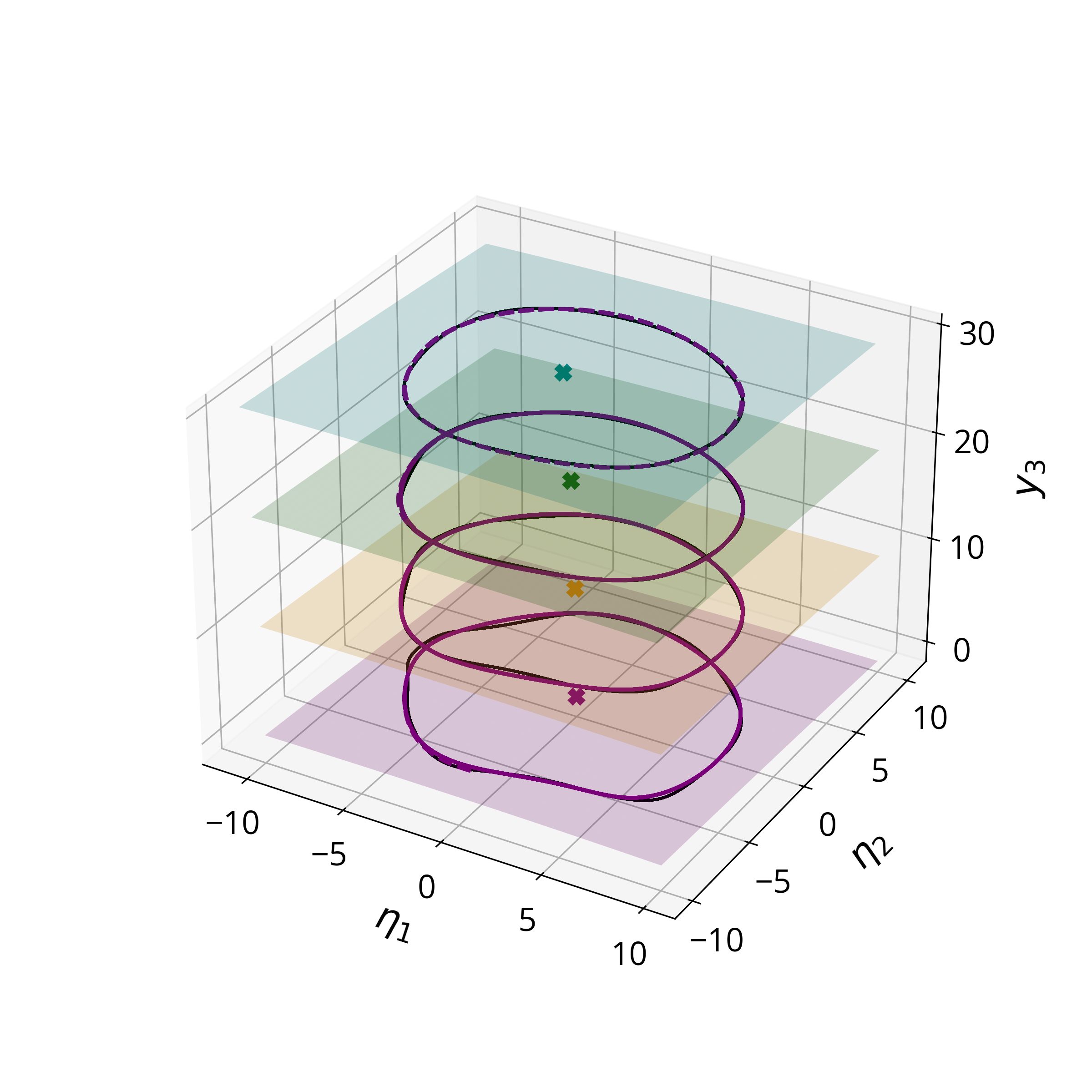}
      \centering
    \includegraphics[width=0.48\linewidth]{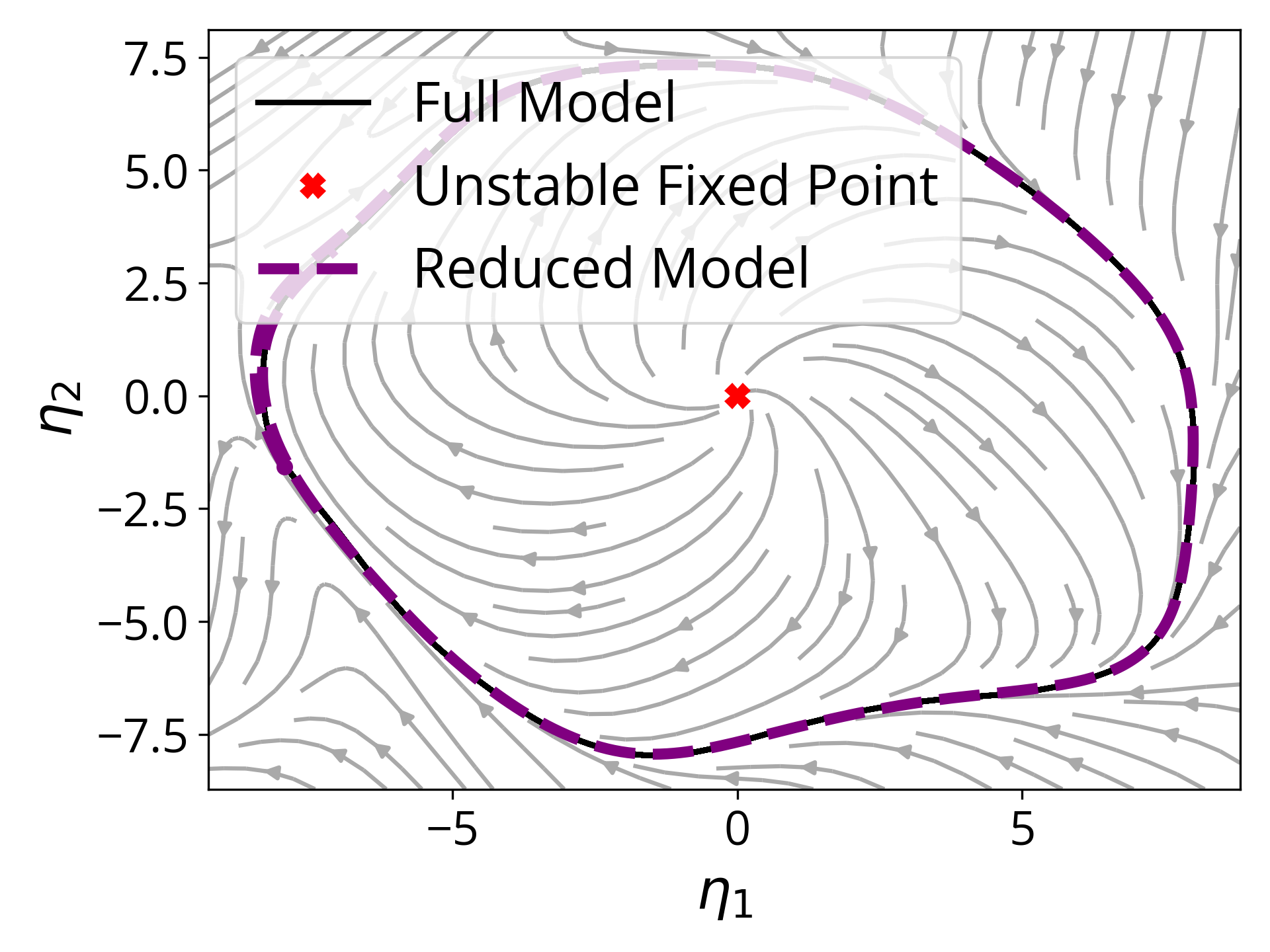}
    \centering
    \includegraphics[width=0.48\linewidth]{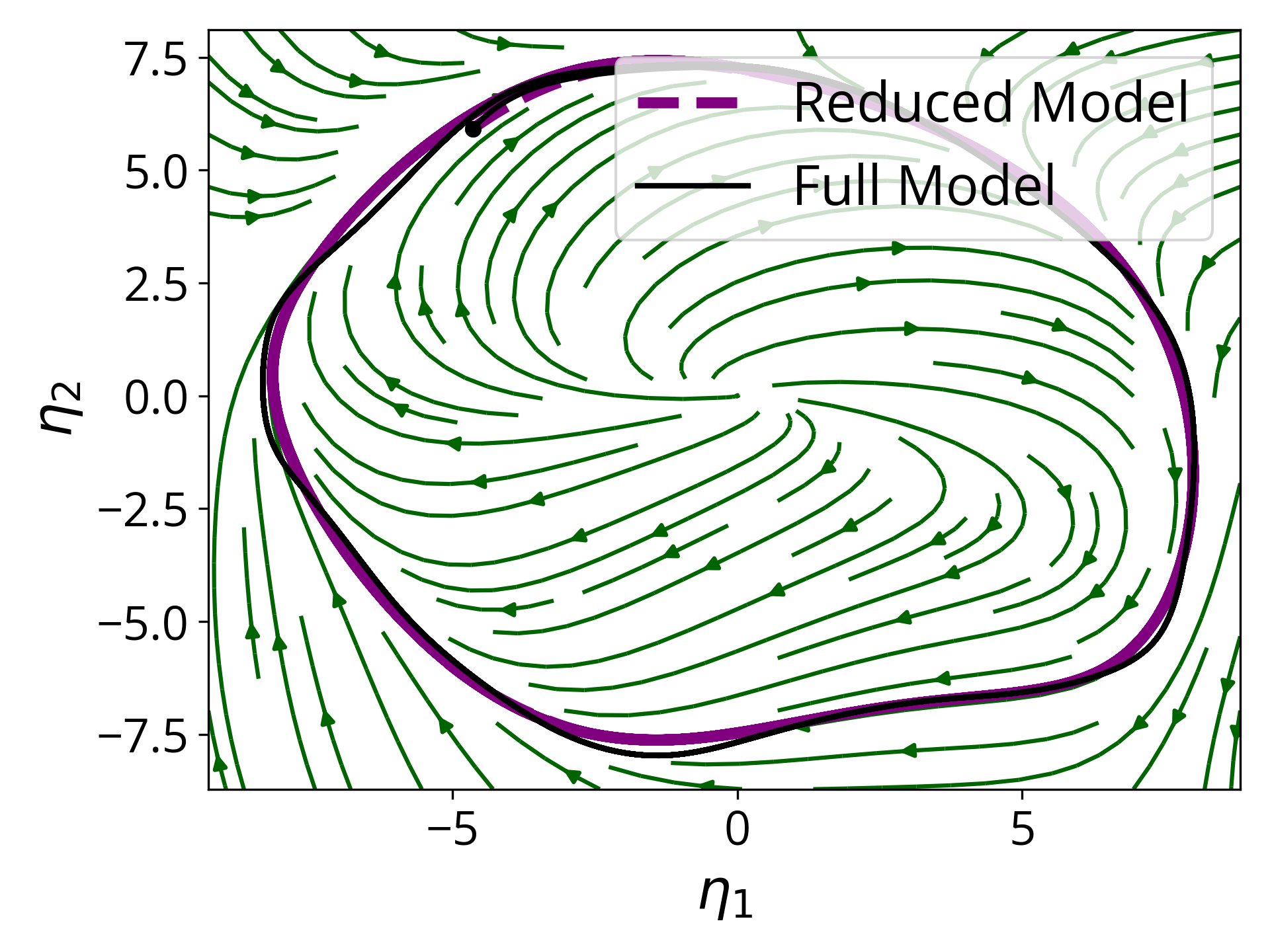}
    
    \caption{SSM carrying the core RNN dynamics and parameter-dependent SSMs with their reduced dynamics projected to the spectral subspaces. (Upper left) Unstable manifold at order three in coordinates $(\eta_1, \eta_2,y_5)$, where $(\eta_1, \eta_2)$ are the coordinates of the unstable subspace. The stable limit cycle on the unstable manifold is globally attracting: in the picture, one full and one reduced trajectory are converging to the stable limit cycle. The MFE and the NMTE calculated on test trajectories have magnitude $\approx 0.02$ and $\approx 0.5$, respectively. (Bottom left ) Projections onto the unstable subspace of the streamlines of the right-hand side of the reduced-order model on the two-dimensional unstable manifold up to order six, with projected full and reduced trajectories.
    (Upper right) Two-dimensional unstable manifolds for different values of the parameter $u$ (and corresponding output frequency values $f$) plotted in the three-dimensional space with coordinates $(\eta_1, \eta_2,y_5)$. The manifolds are shifted on the $y_5$-axis for visualization purposes. (Bottom right) Projection onto the unstable subspace of the streamlines of the right-hand side of the reduced-order model on the two-dimensional unstable manifold up to order six for parameter $u$ corresponding to a frequency of $3.9$ Hz, with projected full and reduced trajectories.}
    \label{fig:swg manifold}
\end{figure}

The unstable manifold, the limit cycle, and full and reduced trajectories are plotted in Fig. \ref{fig:swg manifold} (upper left), while contour lines of the vector field defining the reduced-order model found through SSMLearn at order seven, projected to the unstable subspace $E_2$ can be seen in Fig. \ref{fig:swg manifold} (bottom left). Hence, we obtain a 2D invariant manifold and reduced-order model which provide a complete and detailed description of the dynamics of this high-dimensional, oscillatory RNN. When we choose an input $u$ with frequency $f = 1.9$ Hz, the frequency of the resulting limit cycle captured by the SSM-reduced model is $f_{RED}=1.92$ Hz, which is very close to the measured frequency $f_{FULL} = 1.91$ Hz in the full RNN model.

To analyze the dependence of the limit cycle frequency on the input parameter  $u$, we recall that unstable manifolds depend smoothly on parameters, as long as their underlying unstable spectral subspaces persist under the variation of those parameters (see, for example, \cite{GHDS}). As a consequence, we have a smooth, $u$-dependent family of 2D unstable manifolds in the phase space, as illustrated in Fig. \ref{fig:swg manifold} (upper right).

Along members of this SSM family, we have an explicit, $u$-dependent polynomial expansion for a 2D, parametric, SSM-reduced model of the form

\begin{align}
\begin{split}
\dot{\eta}_1 &= \sum_{i+j+k \leq 4} a_{ijk}^{(1)} \, \eta_1^i \eta_2^j u^k \ + \  \mathcal{O}(5),\\
\dot{\eta}_2 &= \sum_{i+j+k \leq 4} a_{ijk}^{(2)} \, \eta_1^i \eta_2^j u^k \ + \  \mathcal{O}(5).
\label{eq:rom2Dswgpar}
\end{split}
\end{align} The limit cycle frequencies in the parametrized model (\ref{eq:rom2Dswgpar}) agree with the measured frequencies of the full trajectories up to an error of $10^{-1}$.

This example illustrates how our reduction approach can uncover the specific solution that the training process of an RNN selects from many potential solutions to perform a given task. Rather than viewing RNNs as black box systems that replicate some behavior, SSM-reduced modeling uncovers the low-dimensional, interpretable core dynamics \ref{eq:rom2Dswg} of RNNs. This in turn enables us to identify a very specific dynamical structure, a limit cycle, that is directly responsible for producing the output that the RNN was trained to produce.


\subsection{A Model for a Memory-Pro Task}
As our third example, we now describe the SSM-based reduction of the multitasking RNN in \cite{ringattractor} performing a Memory-Pro Task described in the Method section \ref{ss:mtRNN}. This RNN receives as input a fixation input $u_1$ with values in $\{0,1\}$, and two 2D angular sensory inputs
$u_2 = \cos{\theta_1}, \ u_3 = \sin{\theta_1}$, $u_4 = \cos{\theta_2},$ and $ u_5 = \sin{\theta_2}$.
In the setting considered here, the network has to respond according to the angle indicated by the first stimulus when the fixation input is off ($u_1=0$). This is expressed through a three-dimensional readout vector $\mathbf{z}$: $z_1 = u_1$ (fixation) and $z_2 = \cos{\theta_1}$,$z_3=\sin{\theta_1}$.

The authors of \cite{ringattractor} interpret the network dynamics underlying the memory period within the task through a 1D structure composed of {slow points} \cite{badLinearization, mante_context-dependent_2013, slowpoints} located along a circle, which they call a \textit{ring attractor}.
In this structure, the position that a trajectory reaches during the memory period is interpreted as the persistent pattern that the network maintains over time, encoding the correct stimulus it has to respond to when the fixation input is off. 

For this reason, we focus on the input parameter setting corresponding to the memory and context periods, when fixation and context inputs are on and sensory inputs are off. For this configuration, we seek to clarify what dynamical structure is behind what has been empirically labeled as a ring attractor. Such clarification is necessary, as the attractor topology closest to a ring would be that of a torus. There is, however, no observational evidence for sustained quasiperiodic behavior, which would be the hallmark of a toroidal attractor. 

Our analysis reveals an unstable fixed point with two unstable, real eigenvalues and an asymptotically stable fixed point. Based on the spectrum of the linearized dynamics around the unstable fixed point, we conclude that the 2D SSM embodied by the unstable manifold of the unstable fixed point will carry the observable core dynamics once the transients have died out.

\begin{figure}
    \includegraphics[width=0.55\linewidth]{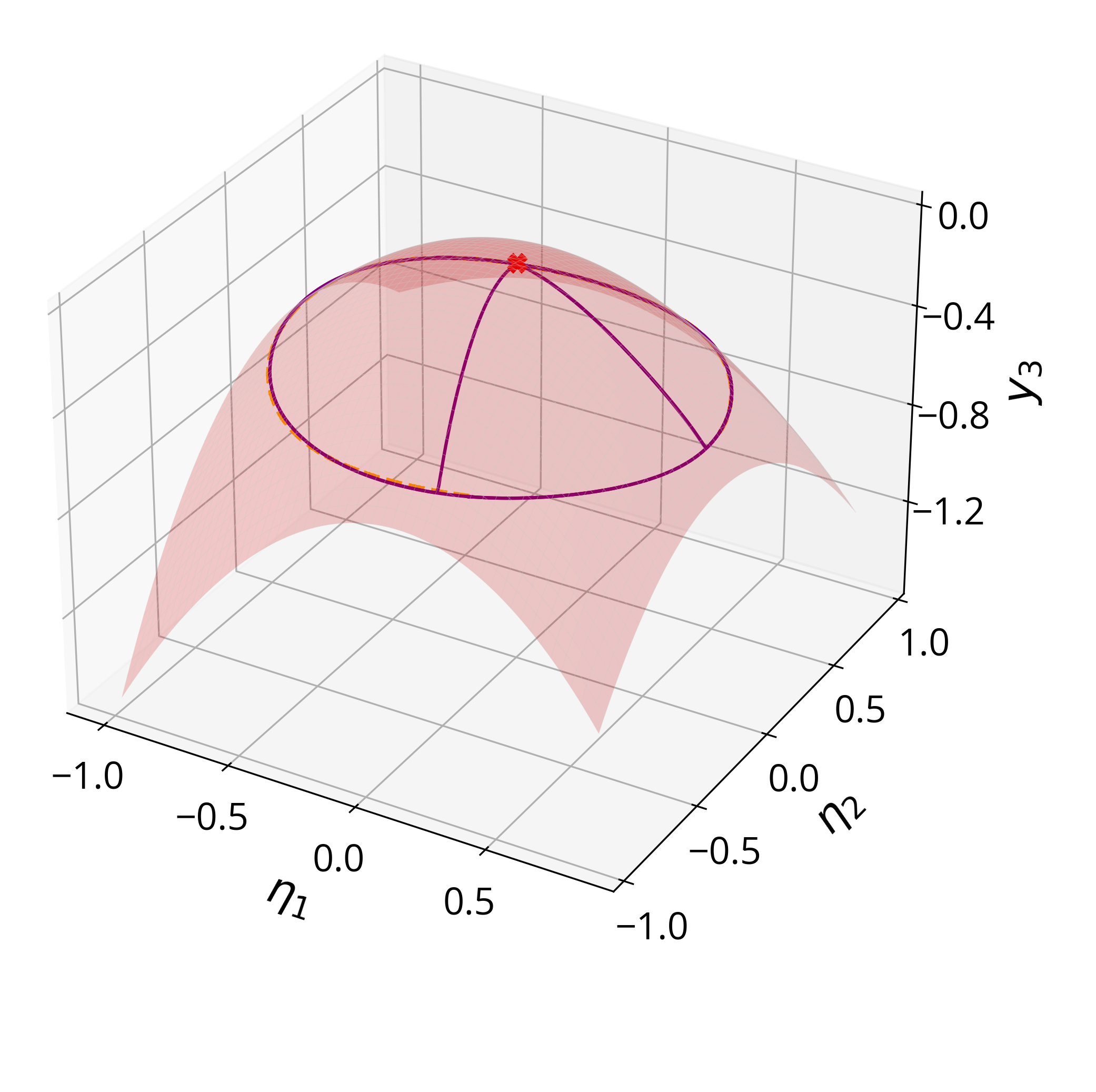}
    \centering
    \includegraphics[width=0.4\linewidth]{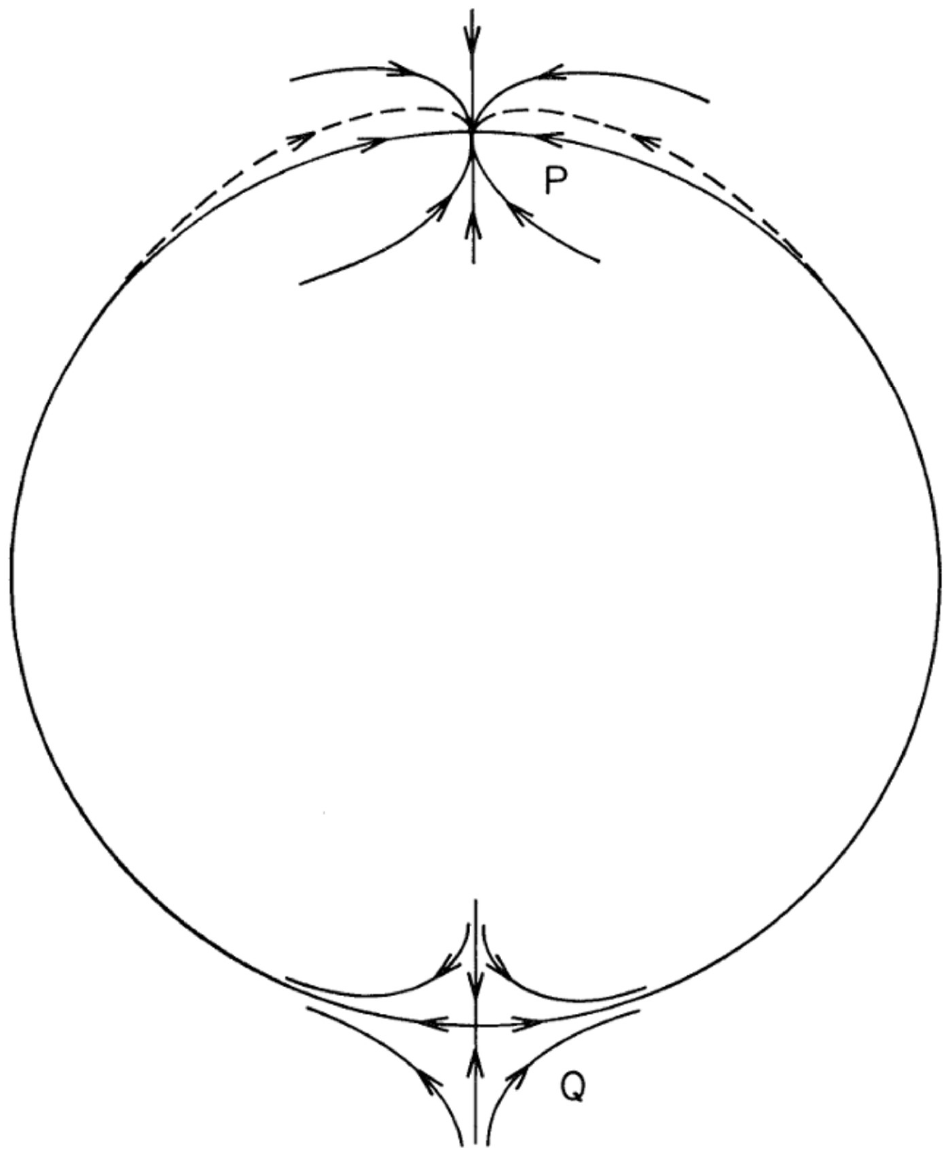}

    \centering
    \includegraphics[width=1\linewidth]{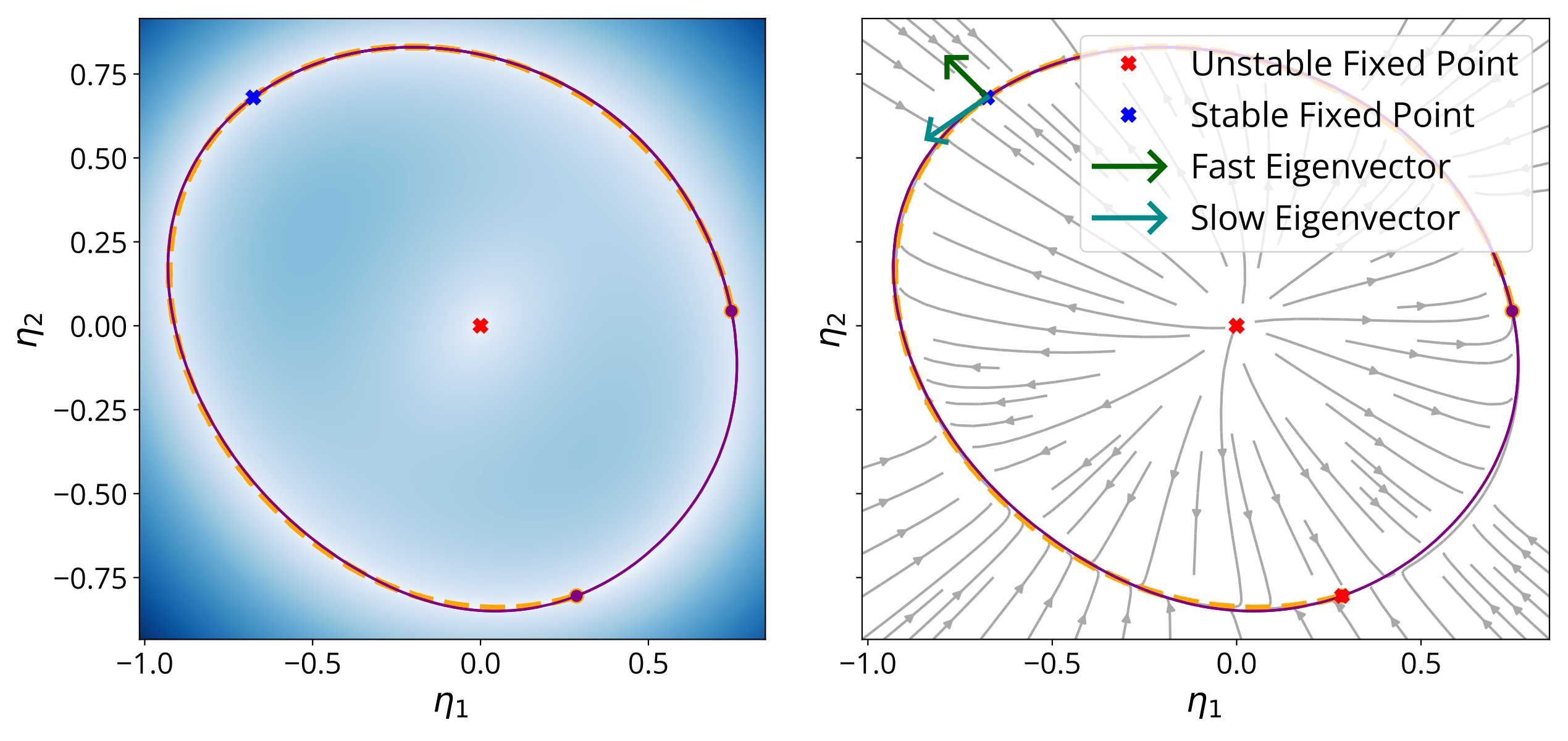}

    \caption{The heteroclinic orbit on the 2D SSM. (Upper left) The two-dimensional unstable manifold attached to the unstable fixed point carries the reduced dynamics. Trajectories (full model trajectories in purple and reduced model in orange) disposed on a circular curve converge to the stable fixed point in infinite time. (Bottom left) Norm of the right-hand side of the reduced model on the unstable manifold plotted on the unstable subspace with coordinates $(\eta_1, \eta_2)$, with darker blues corresponding to larger norm values. The MFE and the NMTE calculated on test trajectories have magnitude $ 10^{-2}$ and $10^{-1}$, respectively.
    (Bottom right) Heteroclinic orbits connecting the unstable fixed point on the bottom right and the stable fixed point on the upper left, projected onto the unstable subspace with coordinates $(\eta_1,\eta_2)$. Streamlines of the right-hand side of the reduced model are projected onto the unstable subspace (in gray). The structure is normally hyperbolic, as we can see from the directions of the slow and the fast eigenvectors of the linearization around the stable fixed point of the full model. (Upper right) Pictorial representation of the normally hyperbolic heteroclinic structure, as an example of normally hyperbolic invariant manifold in \cite{Fenichel}. The unstable manifold of the unstable fixed point Q coincides with the slow stable manifold of the stable fixed point P, and the rate of normal attraction to P is bigger than the rate of tangential compression.}
    \label{fig:slowring}
\end{figure}

We find the unstable manifold using SSMLearn up to polynomial order three (see Fig. \ref{fig:slowring}, upper left). 
 The reduced-order model on this SSM is of the form \begin{equation}
\dot{\eta}_\ell = \sum_{i=0}^3 \sum_{j=0}^{3-i} a_{ij}^{(\ell)} \, \eta_1^i \eta_2^j \ + \  \mathcal{O}(4), \quad \ell = 1,2.
\label{eq:rom2Dmt}
\end{equation}

 Solving for the zeros of the right-hand side of (\ref{eq:rom2Dmt}), we find an additional unstable fixed point. This fixed point is connected to the stable fixed point on the SSM via a pair of heteroclinic orbits that form a structurally stable heteroclinic loop (see Fig. \ref{fig:slowring}, bottom). The structural stability (i.e., smooth persistence under perturbations of the RNN) of this heteroclinic loop follows from the fact that it is normally hyperbolic, i.e., the attraction rates normal to it dominate the attraction rates inside it (\cite{Fenichel}, \cite{Wiggins}).
 
 This strict result from the SSM-reduced model (\ref{eq:rom2Dmt}) clarifies that the true phase space structure behind the empirically observed slow circular structure (ring attractor) is an attracting, 1D invariant manifold. Diffeomorphic to a circle, this manifold is formed by a structurally stable heteroclinic loop between two fixed points inside a 2D slow SSM. 



The structural stability of normally hyperbolic heteroclinic loops explains the ubiquitous presence of slow, circular attractors in RNNs performing working memory tasks (see \cite{nhimsrnn}, \cite{ostojic2023} \cite{memorychad} \cite{Khona2021AttractorAI}). If, in fact, working memory corresponds to the maintenance of persistent neuronal firing patterns, represented in our framework as fixed points in phase space, it is tempting to conclude (see \cite{wm1, wm2, wm3, la1, la2, ringatt_theo, ringatt1, stroud_dyn_coding}) that encoding a continuous variable would require a continuum of such fixed points, each parameterized by a possible value of the variable.

A continuum of fixed points, however, is inherently non-robust and can be destroyed by small perturbations to the neuronal dynamics.
Instead, robust circular attractors that are not limit cycles are necessarily formed by chains of heteroclinic connections among fixed points. Here, we have observed the simplest such chain, one that is formed just by two heteroclinic orbits. 
   
\section{Discussion}

We have shown how analyzing RNNs based on recent results from nonlinear dynamical systems theory reveals the core dynamics of those networks via considerably lower-dimensional, polynomial models on attracting SSMs. Our results firm up previous, mostly heuristic interpretations of RNN dynamics. Specifically, the low dimension of the resulting SSM-reduced models allows us to develop a detailed mathematical understanding of the global RNN behavior. We describe these dynamics in terms of attracting invariant manifolds, the invariant sets they carry, and the transitions among those sets.

In particular, we have obtained a one-dimensional (1D) parameter-dependent model on a \textit{slow manifold} for the 100D, context-dependent decision-making RNN in \cite{mante_context-dependent_2013}. In a similar, data—driven fashion, we have constructed 2D models on slow SSMs for the 100D sine wave generator RNN in \cite{renate}, and the 256D multitasking RNN in \cite{ringattractor} performing a memory-pro task. As opposed to slow manifolds arising in geometric singular perturbation studies (see \cite{fenichelgpt}), SSMs exist without the assumption of a global time scale separation that divides the phase space into slow and fast variables. Therefore, SSMs exist in a more general setting in which the core dynamics settle on a low-dimensional manifold once the faster transient dynamics have died out.

Specifically, for the RNN performing the context-dependent discrimination task, we have explained the “line attractor” observed by \cite{mante_context-dependent_2013} as a robust, parameter-dependent 1D slow manifold. For every configuration of the input parameters, we have identified the corresponding 1D slow manifold (as a section of the parameter-dependent slow manifold in the extended phase space) that captures the reduced dynamics of the full RNN. The resulting behavior following sensory integration arises from the dynamics on a 1D unstable manifold, along which the system converges to the stable fixed point associated with the correct choice.

We have also explained the flexibility of the decision-making process in terms of parameter-dependent SSMs and reduced models. In particular, context and sensory inputs determine which of the two fixed points on the 1D unstable manifold has the larger domain of attraction.

For the sine-wave generator RNN in \cite{renate}, we have uncovered the dynamics responsible for the input-dependent sinusoidal readout with a simple 2D model on a 2D unstable manifold carrying an attracting periodic orbit. This explains how the frequency of this orbit varies with the scalar input through a parameter-dependent model on the parameter-dependent unstable manifold.

Finally, we have found that the empirically described “ring attractor” in the memory-pro task in \cite{ringattractor} is a normally hyperbolic invariant slow manifold that coincides with two heteroclinic orbits connecting a stable and an unstable fixed point on a 2D SSM. The stable fixed point on the heteroclinic circle attracts nearby trajectories, and the heteroclinic structure carries slow, converging dynamics.

The relationship between RNNs trained to perform specific readout functions and the actual computations carried out by neurons is unclear, as many solutions are available for a high-dimensional system to produce a low-dimensional readout. However, an accurate description of the RNNs' phase space geometry enables the generation and verification of hypotheses on the dynamical systems underlying the observed dynamics in the data (see \cite{RNNRev}, \cite{EngelRNN}). Interestingly, in the case of the 100D network trained for a binary choice task (context-dependent decision-making), the solution is a simple 1D bistable system. Similarly, we have found that a limit cycle with input-dependent frequency is responsible for the input-tuned oscillations of the sine-wave generator network. Moreover, the heteroclinic structure explaining the memory-pro task is a robust manifold that appears to be the common mechanism shared by all RNNs trained to reproduce short-term memory mechanisms \cite{nhimsrnn}, \cite{ostojic2023}.

The data—driven SSM—reduction approach used here is general enough to be applied to RNNs performing more complex tasks, whether or not they are related to neural computations. Indeed, RNNs are high-dimensional, nonlinear dynamical systems, whose use extends beyond neuroscience to describe any system with rich dynamics. We rely on simulations of RNN trajectories to find the reduced models of their dynamics through SSMLearn \cite{DataDriven}. This is a fully data-driven algorithm that is applicable to systems whose equations of motion are either unknown or do not display any particular symmetry. Our only assumption is that the RNN has at least one fixed point whose location is at least approximately known. Importantly, the dimension of the SSM—reduced model (typically one or two) is independent of the RNN dimension for generic parameter choices. Therefore,  even for more complex RNNs not considered in this study, the potential for dimensionality reduction using SSMs is virtually unlimited.

The SSM theory was originally developed for model reduction in mechanical engineering, to reduce the dimensionality of high-dimensional mechanical systems, first in an equation-driven fashion \cite{SSMTool}, and later extended to fully data-driven problems \cite{DataDriven}. However, neural systems are inherently different from mechanical systems, which generally have well-understood governing equations that are consistent with observations. For this reason, the application of SSM theory to neural data is promising but not straightforward. Recent extensions of SSM theory cover forced systems \cite{aperiodicSSM} and random systems with small, uniformly bounded noise \cite{randomSSM}. Therefore, a starting point for SSM reduction carried out directly on experimental neural data will require a feasible model for noise in neuronal dynamics.


    \clearpage
\begin{center}
\textbf{\large Methods}

\end{center}
\makeatletter

 
\section{The Theory of Spectral Submanifolds}
   To identify robust, smooth, and unique invariant manifolds in high-dimensional Vanilla RNNs, we employ the theory of SSMs, first formulated by Haller and Ponsioen in \cite{firstSSM} to find invariant manifolds in high-dimensional nonlinear systems around an attracting steady-state solution. The reader can consult \cite{SSMReview} for a complete and detailed introduction to the theory and its applications.

SSM theory is based on a fundamental observation about linear dynamical systems of the form
\begin{equation}
    \dot{x} = A x, \quad x \in \mathbb{R}^N, \quad A \in \mathbb{R}^{N \times N}. \label{eq:linearsyst}
\end{equation}
In this setting, the eigenspaces of \( A \), as well as their algebraic direct sums, are invariant under the flow map and therefore constitute linear and smooth invariant manifolds of the system. We refer to any subspace of the form \( E = E_1 \oplus \cdots \oplus E_k \), where each \( E_j \) is the eigenspace associated with a particular eigenvalue \( \lambda_j \), as a \textit{spectral subspace}. Important examples include the stable, center, and unstable subspaces, spanned respectively by the generalized eigenvectors corresponding to eigenvalues with negative, zero, and positive real parts.

From the eigendecomposition of $A$, we can always find solutions to the linear equations (\ref{eq:linearsyst}) and gather them together to construct infinitely many invariant manifolds tangent to each spectral subspace at the fixed point at the origin. Among these, the spectral subspaces are the smoothest (analytic), whereas all others have limited smoothness governed by the ratios of the real parts of the eigenvalues of $A$.

SSM theory exploits the structure of linear systems and uses linearization theorems to infer the existence, smoothness, and uniqueness of invariant manifolds, SSMs, in nonlinear systems around a fixed point. Although the existence of these manifolds is obtained from linearization, they extend over the domain of linearizability of the system and can carry nonlinear reduced dynamics. We look for SSMs as smooth invariant graphs tangent to the spectral subspaces of the linearized dynamics which, in the general case, are not invariant manifolds of the system themselves and hence not suited for model reduction (see \cite{datadrivenLin} for further discussion).


In particular, consider the nonlinear autonomous system:
\begin{equation}
    \label{eq:DSSimple}
\dot{x} = A x + f_0(x), \quad x \in U \subset \mathbb{R}^N,
\end{equation} where \( f_0 \in C^r \) for some \( r \in \mathbb{N} \cup \{\infty, a\} \) (where $C^a$ indicates that the function is analytic) and satisfies \( f_0(x) = \mathcal{O}(\|x\|^2) \) as \( x \to 0 \). The fixed point at the origin is assumed to be hyperbolic and attracting, meaning that all eigenvalues of \( A \) have strictly negative real parts.

Given a spectral subspace \( E \) for the linearized dynamics, we define its \emph{spectral quotient} as the integer ratio
\begin{equation}
\sigma(E) = \mathrm{Int}\left[ \frac{\max_{\lambda_k \notin \mathrm{spect}(A|_E)} \mathrm{Re}\,\lambda_k}{\min_{\lambda_i \in \mathrm{spect}(A|_E)} \mathrm{Re}\,\lambda_i} \right],
\end{equation}
which compares the fastest decay rate outside \( E \) to the slowest decay rate within \( E \). We can then define a \textit{external nonresonance condition} for E (see \cite{firstSSM}): no low-order ($2\leq m \leq sigma(E)$) integer linear combination of eigenvalues within \( E \) should coincide with an eigenvalue outside of \( E \).

Under the above assumptions, for a nonresonant spectral subspace $E$ of the linearized dynamics of system (\ref{eq:DSSimple}) there exists an invariant manifold, $\mathcal{W}_E(0)$, tangent to \( E \) at the origin. This manifold, the SSM associated with \( E \), is the smoothest among all invariant manifolds tangent to \( E \), and unique in the smoothest class $C^{\sigma(E)+1}$. Furthermore, $\mathcal{W}_E(0)$ inherits smooth dependence on coordinates and parameters from the right-hand side of (\ref{eq:DSSimple}).

In cases where the fixed point is unstable but has a nonempty stable manifold, a similar result holds, provided stronger nonresonance conditions are satisfied, specifically, that no eigenvalue of \( A \) is a nonnegative integer linear combination (with $m\geq2$) of all eigenvalues. Under these assumptions, one still obtains an invariant SSM, $\mathcal{W}_E(0)$, that is unique in its smoothness class and tangent to the spectral subspace, as shown in \cite{haller_fractional}.

In summary, attached to a hyperbolic fixed point with a nonresonant spectral subspace of the linearized dynamics, there is a unique invariant manifold, the SSM, distinguished by its highest smoothness. Since the SSM $\mathcal{W}_E(0)$ is tangent to its underlying spectral subspace $E$ at the fixed point, it is convenient to locally represent it as a graph over that spectral subspace. 

As an example, let $E = \mathrm{span}\{ E_{j_1}, E_{j_2}\}$ be a 2D spectral subspace for $A$ in (\ref{eq:DSSimple}), with $f_0 \in C^a$ (analytic), for which the external nonresonance conditions are met. We assume that we have chosen a basis such that $E$ is parameterized by the reduced coordinates $(\eta_1, \eta_2)$, and the full coordinate system is given by $x = (\eta_1, \eta_2, y_1, \dots, y_{N-2})$. We look for the Taylor expansion of the SSM, $\mathcal{W}_E(0)$, expressed as a graph over $E$ 
\begin{align}
    \begin{pmatrix}
        y_1 \\ \vdots \\ y_{N-2} 
    \end{pmatrix} = h(\eta_1,\eta_2)
    = \begin{pmatrix}
        h_1(\eta_1, \eta_2) \\ \vdots \\ h_{N-2} (\eta_1, \eta_2)
    \end{pmatrix}= \begin{pmatrix}
      \sum_{d=2}^{\infty} \sum_{j=0}^{d} a_{1,d,j} \eta_1^{d-j} \eta_2^{j} \\ 
      \vdots \\
      \sum_{d=2}^{\infty} \sum_{j=0}^{d} a_{N-2,d,j} \eta_1^{d-j} \eta_2^{j}
    \end{pmatrix}
    \label{eq:SSMex}
\end{align}
Since $\mathcal{W}_E(0)$ is invariant, we can write the \textit{invariance PDE} \begin{align} \begin{split}  \dot{x}\Big|_{\mathcal{W}_E(0)} &= A \ \ \begin{pmatrix}
    \eta_1 \\ \eta_2 \\ {h}_1(\eta_1, \eta_2) \\ \vdots \\{h}_{N-2}(\eta_1, \eta_2) 
\end{pmatrix} + f_0(\begin{pmatrix}
    \eta_1&\eta_2&{h}(\eta_1, \eta_2)
\end{pmatrix})  \\&= \begin{pmatrix}
    1 & 0 \\ 0 & 1 \\ \partial_{\eta_1}h_1(\eta_1, \eta_2) &\partial_{\eta_2}h_1(\eta_1, \eta_2) \\ \vdots & \vdots \\ \partial_{\eta_1}h_{N-2}(\eta_1, \eta_2) & \partial_{\eta_2}h_{N-2}(\eta_1, \eta_2)
\end{pmatrix}\begin{pmatrix}
    \dot{\eta}_1 \\ \dot{\eta_2}
\end{pmatrix} \label{eq:invpde}\end{split} \end{align} with boundary conditions at zero ${h}(0)= Dh(0)=0$. Substitution of the Taylor expansion (\ref{eq:SSMex}) into (\ref{eq:invpde}) yields a family of linear equations for the Taylor coefficients $a_{i,j}$, which can be solved recursively under the asumed nonresonance conditions.

 If the right-hand-side of (\ref{eq:DSSimple}) is not analytic but $C^{r}$, $r \in \mathbb{N}\cup \{\infty\}$, increasing the order of the expansion will not necessarily improve the accuracy of the approximation, even on smaller domains. If $f_0$ is $C^r$ with $r$ finite, its expansion would not exist to orders higher than $r$. One should choose accurately the order to truncate the expansion, taking care not to exceed the order at which the approximation ceases to improve on the desired domain.  


This theoretical basis for explicitly constructing invariant manifolds
over arbitrary non-resonant spectral subspaces of the linearized system can be conveniently used for the reduction of high-dimensional nonlinear systems like RNNs. In particular, as a model for the slow, dominant dynamics, we can choose a \textit{slow spectral subspace} $E$ such that $E =  \{E_1 \oplus E_2 \oplus \dots \oplus E_q\}$ is the direct sum of the $q$ slowest eigenspaces whose corresponding eigenvalues have the largest real parts (they represent the lowest decaying dynamics in the linearized system if the fixed point is stable but also growing dynamics if the fixed point is unstable). 

We point out that the dimensional parameter $q$ in the choice of the slow spectral subspace $E$ is a free parameter. It can be chosen so that the slow SSM becomes the low-dimensional attractor that prevails in the available data. For this purpose, we can examine the largest spectral gaps (difference in the real parts) between consecutive (ordered) eigenvalues and pick a slow spectral subspace that is separated from the others on a time scale suitable for the system under study, or simply the slow subspace corresponding to the largest gap.

Back to our example, once a functional expression for the SSM is obtained, we can find the reduced dynamics on the SSM by substituting the $y$ variables with $y|_{\mathcal{W}_E(0)}=h(\eta_1,\eta_2)$ \begin{align}
    \begin{pmatrix}
        \dot{\eta_1} \\ \dot{\eta_2} \end{pmatrix} = \begin{pmatrix}
        f_{0,1}(\eta_1,\eta_2, h_1(\eta_1,\eta_2), \dots,  h_{N-2}(\eta_1,\eta_2)) \\
     f_{0,2}(\eta_1,\eta_2, h_1(\eta_1,\eta_2), \dots,  h_{N-2}(\eta_1,\eta_2))
     \end{pmatrix} \label{eq:reducedeqs}
\end{align}
thus obtaining an explicit polynomial expansion for the 2D reduced-order model on the manifold. Since we have enslaved all the other dynamical variables to $\eta_1$ and $\eta_2$ through the SSM equation (\ref{eq:SSMex}), we can retrieve the dynamics of the full system by inserting the solutions of the reduced equations (\ref{eq:reducedeqs}) into the SSM equation. 

The latest version of the algorithm that performs these equation-driven calculations is available in the open source MATLAB live script package SSMTool \url{https://github.com/haller-group/SSMTool-2.4} (see \cite{jain}).

\subsection{Parameter-dependent and time-dependent models on SSMs}
To understand the input-output relationship in Vanilla RNNs that resemble the sensory input-behavioral output connections in the animal tasks, we must consider how parameter changes affect the dynamics of the RNNs and their SSMs.

In the absence of noise or other types of explicit time dependence (which we may add later as perturbations), the input vector collects the parameters on which the system depends (see the discussion in section \ref{ss:cddmRNN}). Although altering sensory inputs and context inputs (see \cite{mante_context-dependent_2013, badLinearization, pagan}) may affect the geometry of the phase space differently, both inputs keep the RNN to be a autonomous dynamical system.

Specifically, varying the additive parameters within the input vector can modify the phase space geometry by moving the fixed point locations and changing the size of their domains of attraction. Additionally, in cases where the parameter values are close to \textit{bifurcation values} (see \cite{GHDS}) which can lead to changes in the fixed points' stability or to their disappearance. 

If the dependence on the parameters of the right-hand side of the equations of motion is smooth, their SSMs will also display smooth dependency on parameters \cite{firstSSM}, \cite{SSMReview}. However, whenever a fixed point disappears due to a bifurcation, we can no longer rely on the existence of its SSMs as invariant manifolds for the system.

In our examples, fixed points generally do not undergo bifurcations that lead to their disappearance, except in the context of the decision-making RNN when the relevant sensory inputs are varied. In that case, though, the reduced dynamics are carried by a Normally Hyperbolic Invariant Manifold (NHIM) \cite{Fenichel} \cite{Wiggins}. Therefore, the NHIMs theory and the asymptotic properties, uniqueness, and smoothness of unstable manifolds guarantee the persistence of a 1D invariant slow manifold for the whole parameter range considered, even if the fixed points undergo a cascade of saddle-node bifurcations (for details, see the Supplementary Methods \ref{appA:paramSSM}).

In all the other cases, we can rely on the smooth dependence of SSMs on parameters to find parameter-dependent reduced-order models for the RNN dynamics. In fact, we can extend the phase space of parameter-dependent dynamical systems $\dot{x} = f(x; \mu)$, with $\mu \in \mathbb{R}^p$ the parameter vector, by treating the parameters as \textit{dummy variables} for the system, which have trivial dynamics $\dot{\mu}=0$. We can therefore find polynomial parameter-dependent expansions for SSMs and reduced-order models by solving the invariance equations in the extended phase space, including the parameter directions in the spectral subspaces. 

On the other hand, to examine the impact of bounded noise on the dynamics of recurrent neural networks (RNNs), we apply an extension of SSM theory to nonautonomous systems subjected to weak aperiodic forcing, as discussed in \cite{aperiodicSSM}.  We treat uniformly bounded noise as a form of discontinuous-in-time, weak forcing and derive time-dependent SSMs and reduced-order models that capture the noisy dynamics (refer to Supplementary Methods \ref{appB} for the general discussion on non-autonomous SSMs). These models will maintain the same phase space geometry as the autonomous system at any given moment, but they will evolve over time due to the influence of noisy perturbations.

In particular, fixed points for the autonomous system perturb, in the noisy setting, into noisy \textit{anchor trajectories} with the same stability properties.
 We can think of such a trajectory as the perturbation of the line that would be present in the extended phase space with coordinates ${X}= ({x},t)\in\mathbb{R}^{N+1}$ with dynamics ${\dot{X}}=({\dot{x}},1), \ L = \{({x},t) \in\mathbb{R}^{N+1}\ | \ {x}= {0},\ \ t\geq0\}$.

 Any autonomous SSM perturbs into a time-dependent SSM that evolves over time along the anchor trajectory perturbing from the fixed point the autonomous SSM was attached to. The time-dependent SSM can be seen as the result of perturbing with additive, small noise terms the autonomous SSM in the extended phase space $\mathbb{R}^{N+1}$ with coordinates $({y},t)$. A more detailed discussion and the calculations specific for the noisy RNN systems can be found in the supplementary methods \ref{appB:anchor} and \ref{appB:SSMtdep}.

\subsection{Data-Driven Computations of SSMs}
As we have seen, SSMs can be found locally as invariant graphs over their tangent spaces at the fixed point and approximated through Taylor expansions. This makes it possible to use a simple data-driven algorithm to obtain low-dimensional, predictive models for high-dimensional, nonlinear systems for which the equations of motions are not known (see \cite{DataDriven}).The data-driven algorithm that finds SSMs and reduced models on them, SSMLearn, is open-source and available at \url{https://github.com/haller-group/SSMLearn}. 

In this study, we combine the equation and data-driven approaches, since we use simulated RNN trajectories to find optimal polynomial SSMs and reduced-order model approximations through SSMLearn.

More specifically, SSMLearn is an algorithm designed to find the optimal polynomial coefficient matrix $H$ for the SSM equations 
\begin{equation} y = V_EV_E^Ty+H\phi_{d,2:M}(V_E^Ty) ;\label{eq:ddssm}\end{equation}
where $d$ is the dimension of the SSM and $M$ is the order of approximation of the polynomial expansion. The function $\phi_{d,2:M}(q)$ outputs the vector of all $m_M$ monomials from order two up to order $M$ given a $d$-vector $q$, $H$ is the $N\times m_M$ coefficient matrix, $V_E^T$ is the projection matrix to the spectral subspace,  $\eta: = V_e^T \ y$ is the $d$-dimensional coordinate vector projected to the spectral subspace.
 
The optimization problem for $H$, and possibly the $N\times d$-matrix $V_e$, is formulated as minimizing a quadratic cost function that measures the distance between train trajectories and the SSM. When an estimate of the spectral subspace is already available ($V_e$ is known), the optimization problem is convex.

It is then possible to compute the polynomial expansion for the reduced equations of motions $\dot{\eta}$, $\eta = V_e^Ty$. In particular, SSMLearn finds the optimal polynomial coefficient matrix $W_r$ for the equation
\begin{equation} \dot{\eta} =  W_r\phi(\eta)_{d,1:M_r} \end{equation}  minimizing the distance between trajectories $\eta(t)$ and the train trajectories projected to the spectral subspace. Finally, full-order trajectories can be obtained by lifting the reduced-order trajectories through the SSM equations.

 \section{The Vanilla RNN Models}
\label{ssec::RNNDS}
The RNNs described in \cite{renate} and \cite{ringattractor} are trained to obtain specific readout dynamics associated with behavioral experiments. These RNNs are designed to produce the correct output values based on a given set of input parameters, enabling them to replicate the target behaviors exhibited by trained animals. In this context, the optimal RNN parameters are identified to accurately reproduce a $N_o$-dimensional readout function $\mathbf{z}$ that depends on the specific task.

The network units are coupled through a connectivity matrix $\mathbf{W}$ and have nonlinear dynamics due to a bounded, smooth, nonlinear activation function. They have equations
\begin{align}
    \begin{split}
        \label{eq:VanRNN}
        \tau\, &\dot{\mathbf{x}}=-\mathbf{x} + \mathbf{W}\,\mathbf{r(x)} +  \mathbf{B}\,\mathbf{u}(t) +\mathbf{\sigma}(t), \\ \mathbf{x} &\in \mathbb{R}^N, \ \ N = 100\ , \ \  \mathbf{r(x)}= (\tanh(x_1),\dots \tanh(x_N))^T \in \mathbb{R}^N,  \\ \mathbf{W} &\in \mathbb{R}^{N\times N}, \,\,\mathbf{B} \in \mathbb{R}^{N\times N_i},\,\, \mathbf{u} \in \mathbb{R}^{N_i}
    \end{split}
\end{align}
where $\mathbf{u}$ is the input vector of dimension $N_i$ and  $\mathbf{\sigma} \approx\, A\,\,\mathcal{N}(0,1) \in \mathbb{R}^{N}$ is a vector with entries drawn at each time step from a Normal distribution with zero mean and unit standard deviation\footnote{In our treatment, we will reject drawings outside $\pm  3$ standard deviations from the mean to make the Gaussian noise bounded.}.

The matrices $\mathbf{W}, \ \mathbf{B}$ and $\mathbf{Y} \ \in \mathbb{R}^{N_o\times N}$ ($N_o$ is the readout vector dimension) are selected by training to produce a \textit{ read-out function} $z(t) = \langle \mathbf{Y}\, ,\,\mathbf{r}(\mathbf{x}(t)) \rangle$ with specific dynamics adapted to each task.

We first remove the noisy components from equations (\ref{eq:VanRNN}) and choose constant input vectors (which, in the experiments in \cite{mante_context-dependent_2013}, \cite{renate}, \cite{ringattractor}, are only varied among trials). We then vary the input parameters to obtain manifolds and reduced ODEs that explicitly depend on them and describe the parametric changes to the geometry of the phase space. Finally, we insert noise back into the equation to obtain a noisy, time-dependent model. 

\subsection{RNN for Context-Dependent Decision-Making}
\label{ss:cddmRNN}
We now describe the specific features of the \textit{Vanilla RNN} described in Eq. \ref{eq:VanRNN} that adapt it to the context-dependent decision-making task of Mante et al. in \cite{mante_context-dependent_2013}.

In the experiment reported in \cite{mante_context-dependent_2013}, monkeys were trained to distinguish between the direction of motion (left or right) or the color (red or green) of a random-dot display based on the given context cue. Meanwhile, extracellular responses from the monkey's neurons were recorded in and around the frontal eye field and analyzed in the space of the first three principal components. 

\begin{figure}
    \centering
    \includegraphics[width=0.45\linewidth]{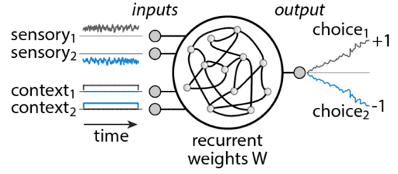}
    \centering
    \includegraphics[width = 0.4\linewidth]{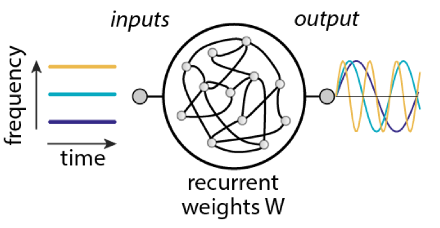}
    \caption{Pictorial representation of the context-dependent task performed by the RNN in \cite{mante_context-dependent_2013} and \cite{renate} (right). Pictorial representation of the sine-wave generator RNN in \cite{renate} (left).}
     \label{fig:cddRNN}
\end{figure}

At each trial, neurons in the RNN receive two independent sensory inputs, each with a positive or negative value, that mimic evidence for motion and color in a single random dot stimulus, and the context inputs (see Fig. \ref{fig:cddRNN}, right) indicate which of the two sensory inputs will determine the sign of the final scalar readout $z \in \mathbb{R}$ ($N_o = 1$).  

The $100$D network evolves according to Eq. (\ref{eq:VanRNN}), with a 4D ($N_i = 4$) input vector $\mathbf{u} = (
    u_1 \ u_2 \ u_3 \ u_4)^T$ that accounts for both context and sensory stimuli in the experiment. In particular, its first two entries $u_1 = s_1$ and $u_2 = s_2$ ( represent sensory stimuli (color and direction of motion, respectively). The network was trained on the set of values $\{-0.15,-0.036,-0.009,0.009,\\0.036,0.15\}$. The last two entries in the input vector $u_3$ and $u_4$ are one-hot-encoded and represent context 1 when $u_3 = c_1 = 1 $ and $u_4 = c_2 = 0$, and context 2 when $u_4 = c_2  =1$ and $u_3=c_1 =0$. 
    
    The input vector $\mathbf{u}$ is piecewise constant during one trial: in the first part of the trial, before the \textit{burnlength time} $T_B$, the sensory inputs are switched off $s_1 =s_2 =0$ and only one of the context inputs is on. During this period, the target readout $z^*(t) = const$.  Instead, when $t >T_B $, the sensory inputs are switched on $u_1 = s_1$, $u_2 =s_2$ and the sensory integration process begins. The target readout is a scalar function converging to either $+1$ or $-1$, depending on the selected sensory input.
    
    We will treat this discontinuous-in-time dynamical system as the union of two separate systems, one for each time interval.
   
\subsection{Sine-Wave Generator RNN}
\label{ss:swgRNN}
The sine-wave generator network is a $100$D Vanilla RNN with equations of motions of the same structure as in Eq. (\ref{eq:VanRNN}). The behavior implemented in this case is not linked to a specific experiment, but is rather a general model for oscillatory responses, which are of interest in the study of muscle contraction (see \cite{renate}).

The input $u$ is, in this case, a scalar ($N_i =1)$, whose values represent the desired oscillation frequencies. The network is trained so that the scalar readout $z(t) = \langle \mathbf{Y}, \mathbf{r}(\mathbf{x}(t)\rangle$ evolves in time as a sine-wave with the frequency indicated by the scalar input (see Fig. \ref{fig:cddRNN}, left).

\subsection{Multitasking RNN: Memory-Pro Task}
\label{ss:mtRNN}
The multitasking RNN implemented in \cite{ringattractor} is \begin{equation}
    \label{eq:MultitaskRNN} \tau \dot{\mathbf{x}} = - \mathbf{x} + \mathbf{r}(\mathbf{W}_{rec}\mathbf{x}+ \mathbf{W}_{in}\mathbf{u} + \mathbf{B}_{rec}),
\end{equation} a variant of the Vanilla RNNs described in Eq. \ref{eq:VanRNN} with $\mathbf{x} \in \mathbb{R}^{256}$. The matrices $\mathbf{W}_{rec}, \ \mathbf{W}_{in}, \ \mathbf{B}_{rec} $ are trained so that the RNN can perform multiple tasks, depending on the values of the 20D input vector $\mathbf{u}$. The one-hot-encoded last fifteen entries $u_i, \ i = 6, \dots, 20$ indicate which task the network has to perform in each trial. The first entry of the input vector plays the role of a fixation input, taking values 0 or 1 depending on whether the network, similarly to animals in an experiment, has to respond to sensory inputs or not. Finally, the two 2D angular sensory inputs are represented in the entries
$u_2 = \cos{\theta_1}, \ u_3 = \sin{\theta_1}$, $u_4 = \cos{\theta_2},$ and $ u_5 = \sin{\theta_2}$.

When the fixation input is off ($u_1=0$), the network has to respond according to the angle indicated by the stimulus that is relevant to the current task. The three-dimensional readout vector $\mathbf{z}$ ($N_o = 3$) encodes whether the network is fixating or not in its first entry $z_1 = 1$ or $z_1 =0$, and for the output angle through the values of $z_2=\cos\theta_{out}$ and $z_3 = \sin\theta_{out}$. 

\begin{figure}
    \centering
    \includegraphics[width=0.4\linewidth]{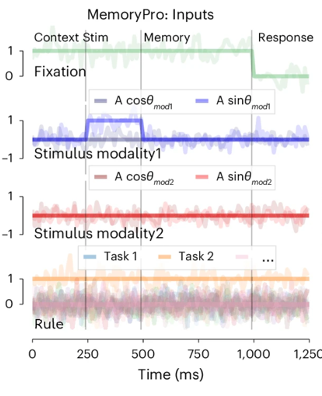}
    \includegraphics[width = 0.45\linewidth]{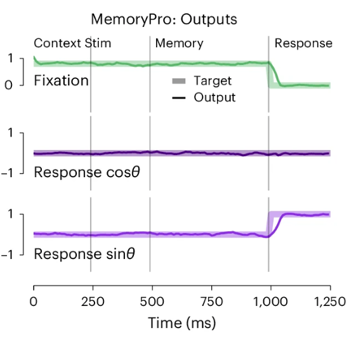}
    \caption{Pictorial representation of the inputs (right) and outputs (left) to the multitasking RNN in \cite{ringattractor} performing a Memory-Pro task.}
    \label{fig:multitaskRNN}
\end{figure}

In a Memory-Pro task, the RNN is exposed to piecewise constant changes in the input vector. In the context period $[t_0, t_s]$, the sensory inputs are off, while the fixation and context inputs are on. Then, in the stimulus period $t_s, t_m$, sensory inputs are switched on, while all the other inputs remain constant. Sensory inputs are switched off again in the memory period $[t_m, t_r]$, after which the fixation stimulus is switched off, and the network has to output the correct angle $\theta_{out}$ (see Fig. \ref{fig:multitaskRNN}). 

This task is designed to understand possible dynamical motifs that underlie working memory in animals. 
In our analysis, we will focus on the configuration during context and memory periods, with context and fixation input on and sensory inputs off. 

    \printbibliography[heading=bibintoc, title={Bibliography}]

\clearpage
\begin{center}
\textbf{\large Supplementary Information}
\end{center}
\setcounter{equation}{0}
\setcounter{figure}{0}
\setcounter{page}{1}
\setcounter{section}{0}
\makeatletter
\renewcommand{\theequation}{S\arabic{equation}}
\renewcommand{\thefigure}{S\arabic{figure}}
\renewcommand{\thesection}{S\arabic{section}}



\section{Context-Dependent Decision-Making RNN: Details on Reduced Dynamics}
\label{appD:const}
We begin our analysis of the dynamical system in Eq. (\ref{eq:VanRNN}) by removing the noise components and fixing the sensory and context inputs. In particular, for the computations shown in this section, we choose $u_1 = s_1 =  u_2 = s_2 = 0$ and $u_3 = c_1  =0$, $u_4= c_2 = 1$ when $t<T_B$ (sensory inputs are off and we choose context 2) and we switch on the sensory inputs to $u_1 =s_1= 0.036 $, $u_2 = s_2 = 0.15$ when $t>T_B$.

 We obtain two autonomous dynamical systems with constant inputs, one for $t< T_B$ and the other for $t>T_B$. When $t<T_B$, we expect to find evidence of the set of \textit{slow points} in the phase space found by Mante et al. in \cite{mante_context-dependent_2013}. When $T>T_B$, we expect this line of slow points to disappear and to uncover some mechanism that can account for the binary choice performed by the network according to the input vector's entries. In both cases, a reduced-order model on a one-dimensional SSM can explain the dynamics of the network. 

For the sensory-inputs-off period $t< T_B$, we find a stable fixed point in the phase space region explored by trajectories starting from initial conditions within a ball $B_{5 \delta}$, $\delta = 0.1$. We are interested in this specific region of the phase space, since the network in \cite{renate} and \cite{mante_context-dependent_2013} was always initialized at zero, and the bounded noise contributions have maximal amplitude $3\delta$. The network was trained to obtain a constant readout dynamics, so convergence to a stable fixed point is a simple realization of this behavior. We also find the two additional stable fixed points corresponding to the two choices, which turn out to be present for all parameter values explored in the task (see Fig. \ref{fig:bif_diag}), but trajectories do not explore their domains of attraction in the pre-stimulus period.

\begin{figure}
    \centering
    \includegraphics[width=0.48\linewidth]{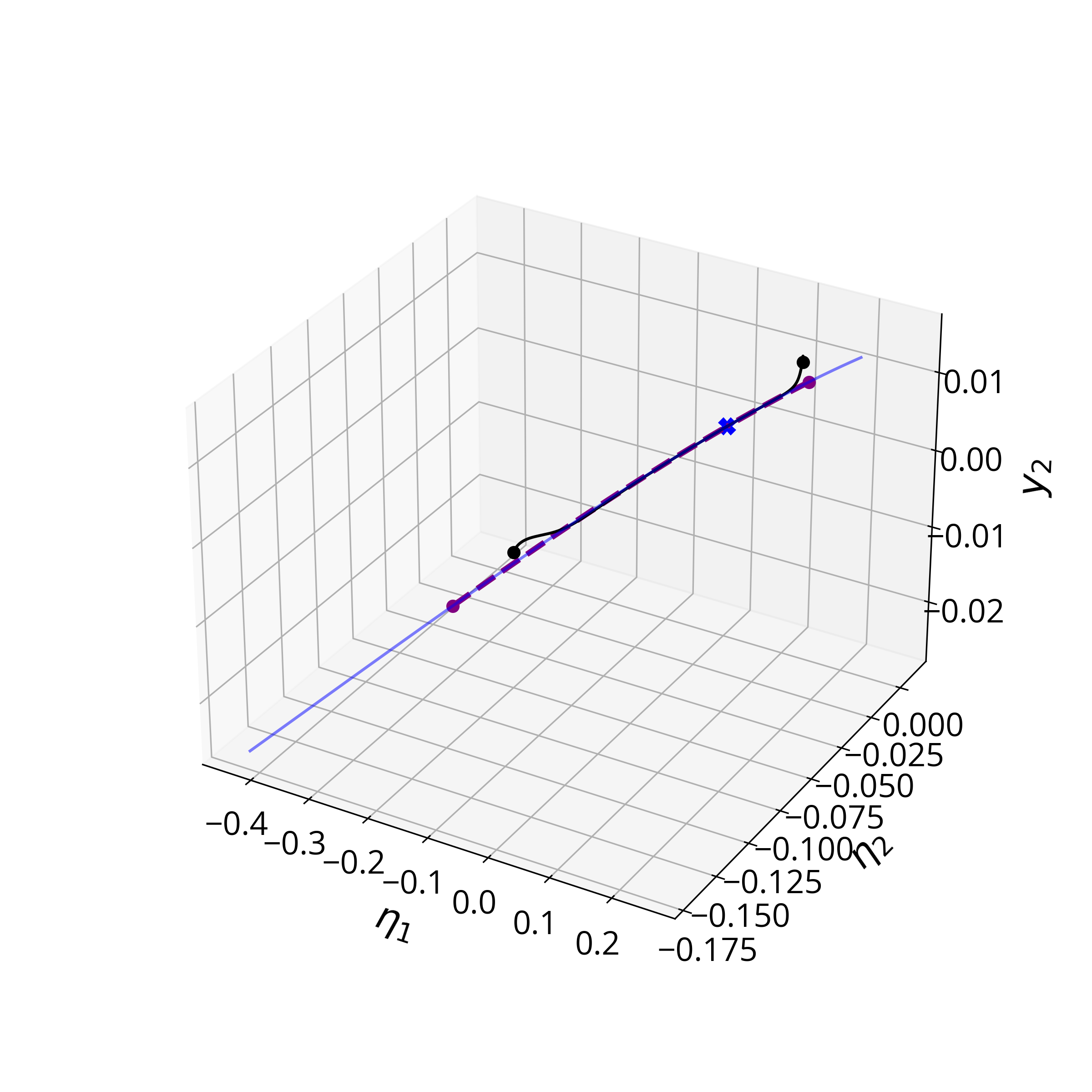}
    \centering
    \includegraphics[width=0.48\linewidth]{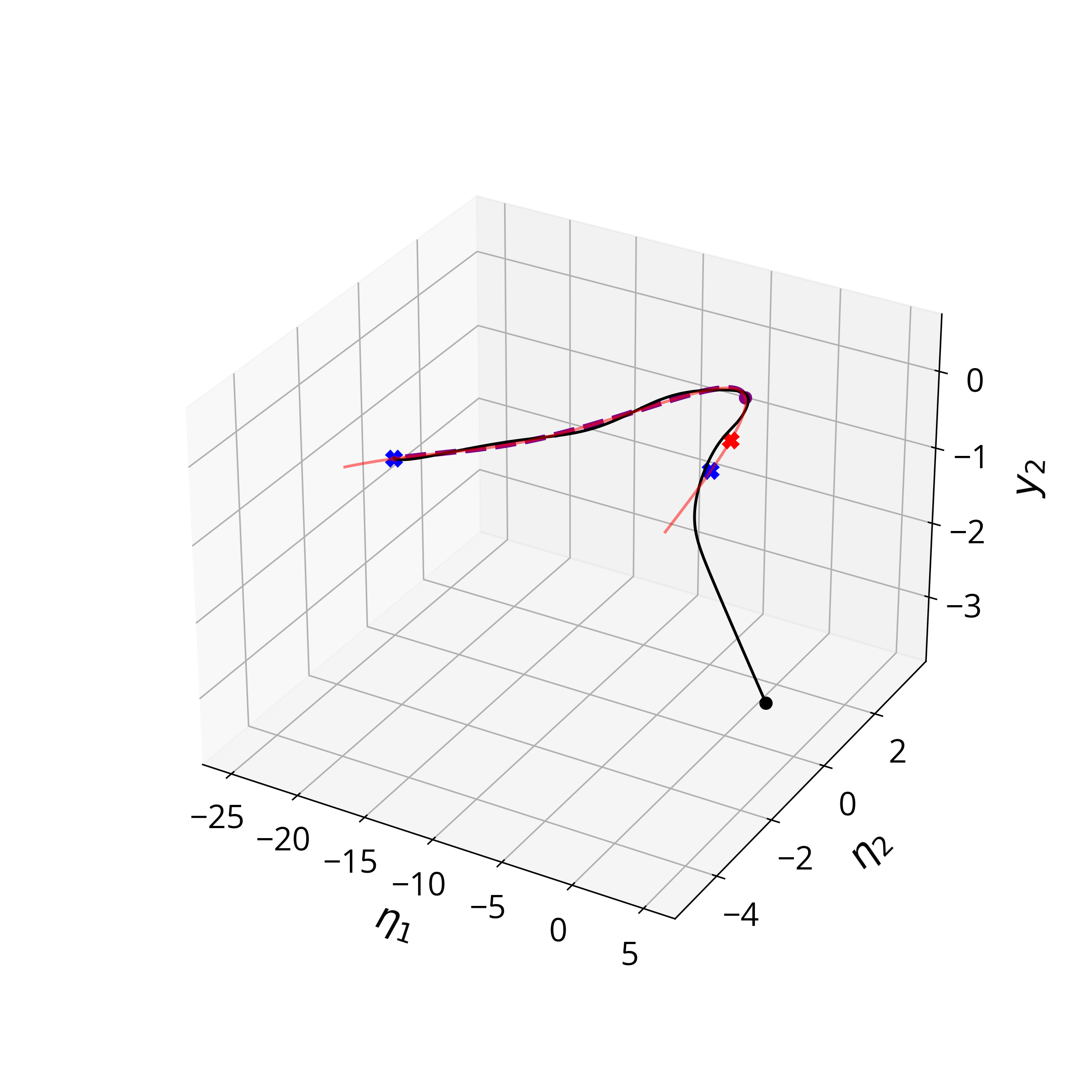}
      \centering
    \includegraphics[width=0.48\linewidth]{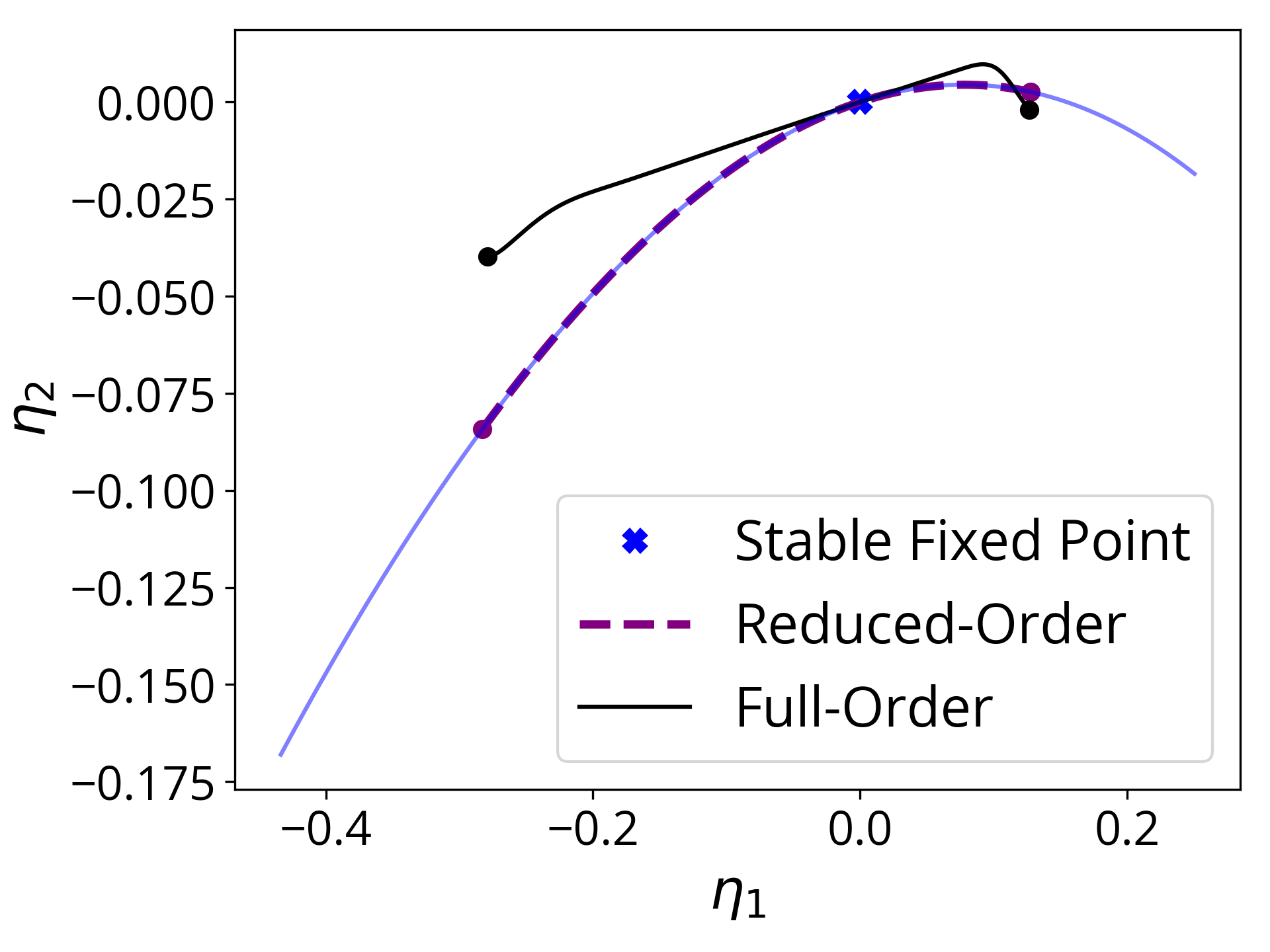}
    \centering
    \includegraphics[width=0.48\linewidth]{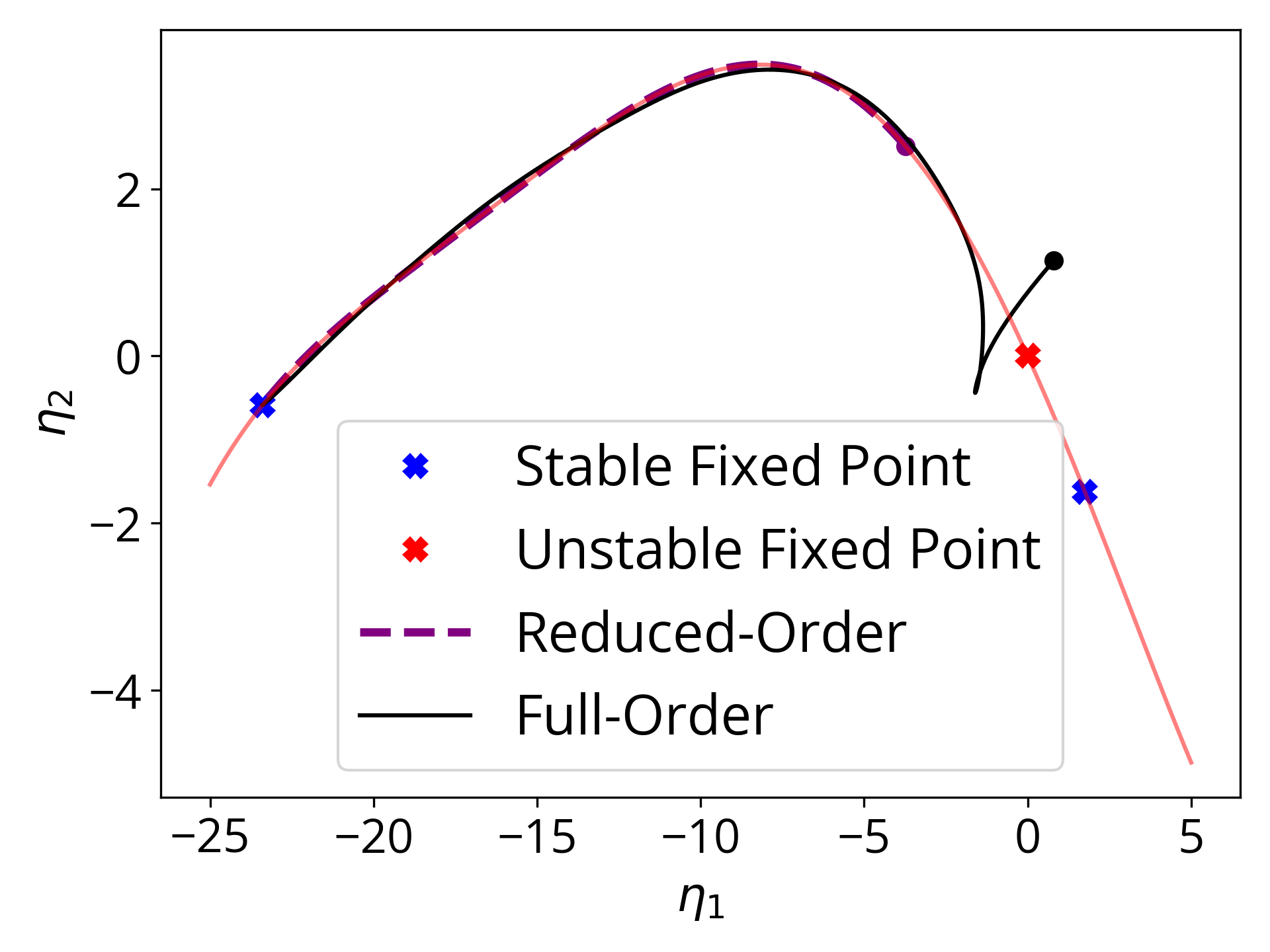}
    \caption{(Top) The slowest SSM at order five (left, sensory inputs off, $s_1 = s_2 =0$, context 2 $c_2 = 1, \ c_1 = 0$) and the unstable manifold of the unstable fixed point at order five (right, sensory inputs on $s_1 = 0.036, s_2 = 0.15$) and test full-order and reduced-order trajectories plotted with coordinates $(\eta_1, \eta_2, y_2)$, where $\eta_1, \eta_2$ parametrize the spectral subspace $E_2$ (corresponding to the first two Principal Components) and the $y_i$, $i = 1, \dots, 100$, are coordinates for the 100D phase space $\mathbb{R}^N$, where the origin is at the anchor point. (Bottom) Manifolds projected on the spectral subspace $E_2$ of the unstable fixed point, spanned by the eigenvectors corresponding to its slowest eigenvalues (one stable and one unstable), with coordinates $(\eta_1, \eta_2)$, corresponding to the first two Principal Components space.}
    \label{fig:SupSSMs}
\end{figure}

The linearized dynamics around the stable fixed point has eigenvalues with a large gap between the real part of the first eigenvalue $\lambda_1$ and all the real parts of the other eigenvalues, which are smaller. Based on this observation regarding the linearized dynamics, we can expect that, at the considered timescales, the dynamics will quickly decay to the (existing, unique and smooth) 1D slowest Spectral Submanifold tangent to $E_1 = \mathrm{span}\{\mathbf{e}_1\}$ (eigenvector of eigenvalue $\lambda_1$). We find the SSM pictured in Fig. \ref{fig:SupSSMs} (left), and Fig. \ref{fig:unstablemfld}, at fifth order and the reduced-order model for the network dynamics on it through SSMLearn. In Fig. \ref{fig:unstablemfld} (center left), we plot the right-hand side of the reduced 1D equations of motions up to the fifth order.

Trajectories quickly converge to the slowest SSM. Once they have approached the manifold, they slowly approach the stable fixed point. This creates the impression of a set of "slow points" moving along an invariant curve in the phase space, converging slowly to a stable fixed point. Therefore, the {line attractor} reported by Mante et al. in \cite{mante_context-dependent_2013}, is just the slowest SSM attached to the stable fixed point: it attracts nearby trajectories and drives slow dynamics due to the (infinite-time) convergence to the stable fixed point. 

We now turn to the analysis of the sensory integration period $T>T_B$ and set $\mathbf{u} = (0.15,0.15, 0,1)^T$. We find one unstable fixed point and two stable fixed points, i.e., two-point attractors in the phase space. Trajectories with random initial conditions are more likely attracted to the fixed point with the larger domain of convergence.

The unstable fixed point has a single unstable eigenvalue, hence a 1D unstable manifold. Since the gap between the real parts of the unstable eigenvalue and stable eigenvalues is wide, we expect the unstable manifold to capture the core dynamics of the system. 

We find the unstable manifold with SSMLearn up to order five (see Fig. \ref{fig:unstablemfld}, upper left, and Fig. \ref{fig:SupSSMs}, right), and the 1D reduced-order model on the unstable manifold (see Fig. \ref{fig:unstablemfld}, center left), through which it is possible to simulate test trajectories and predict the whole system's dynamics.       

To complete our analysis, we add noise to the equations of motion (\ref{eq:VanRNN}). For illustrative purposes, we focus on the sensory input-off time interval when there is a stable fixed point in the explored region of the phase space, and the core dynamics are captured by a reduced-order model on the slowest SSM. We also choose context 2 in the input vector $\mathbf{u}$ to be consistent with the results above for the autonomous network.

We apply the results of \cite{aperiodicSSM} for weak, time-dependent, uniformly bounded, and possibly discontinuous-in-time forcing terms. With this aim, we treat the noise term $\mathbf{\sigma}(t) \approx\, 0.1\,\,\mathcal{N}(0,1) \in \mathbb{R}^{N}$ as small, discontinuous bounded forcing that depends on time and is independent of the position in the phase space.

The stable fixed point $\mathbf{x}_0$ for the autonomous system perturbs into an attracting, uniformly bounded anchor trajectory. We approximate this anchor trajectory through the formal expansions derived in \cite{aperiodicSSM}. We stop at order two, since the resulting approximated trajectory at order one is already a good approximation of the attracting anchor trajectory that we could estimate by evolving an initial condition in $\mathbf{0}$ with the noisy equations of motions \ref{eq:VanRNN}. The formulas for these calculations can be found in the Supplementary Materials \ref{appB:anchor}.

Once we have a good approximation of the anchor trajectory, we can compute the aperiodic, time-dependent slowest Spectral Submanifold $\mathcal{W}(E_1,t)$ attached to $\mathbf{y}^*(t) = \mathbf{x}^*(t)-\mathbf{x}_0$ when the sensory inputs are off and tangent, for every fixed $t$ to the spectral subspace $E_1 =\mathrm{span}\{\mathbf{e}_1\}$ in $\mathbf{y}^*(t)$.
We truncate the approximation for $\mathcal{W}(E_1,t)$ at order $M =2$. The calculations for the first and second-order terms of the expansion can be found in the Supplementary Materials \ref{appB:SSMtdep}. 

The perturbed dynamics consist of fast convergence to the (moving) slowest stable SSM followed by slow convergence to the moving anchor trajectory. 

\section{Context-Dependent Decision-Making RNN: Parameter Dependent SSMs}
\label{appD:par}

We recall that the 4D input vector $u$ in the decision-making RNN of \cite{mante_context-dependent_2013} has entries $u_1 = s_1$, $u_2 = s_2 $ corresponding to sensory inputs, and one-hot-encoded entries $u_3= c_1$, $u_4 = c_2$ playing the role of context inputs. Here, we study changes in the phase space geometry due to parametrically varying both context and sensory inputs, and highlight how the role of the input vector in changing the RNN dynamics can be explained, in this case, in terms of parameter-dependent variations in the phase space.

We will consider a 1D parameter $\mu$, corresponding to either context or relevant sensory/irrelevant sensory inputs, and compute the $\mu$-dependent 1D slow manifolds and reduced-order models by looking for the 2D smooth invariant manifold in the extended phase space $\mathbb{R}^{N+1}$ with coordinates $(\mathbf{y}, \mu)$, tangent to the two-dimensional plane spanned by $\mathbf{e}_1$ (the eigenvector corresponding to eigenvalue $\lambda_1$ as above) and the vector spanning the parameter direction $\mathbf{c}_{N+1}$ in the extended phase space, at the point $(\mathbf{0}, 0)$. 

We seek a such inveriant manifolds as invariant graphs over the plane $P = \mathrm{span}\{\mathbf{e}_1, \mathbf{c}_{N+1}\}$ with coordinates $\eta_1$ and $\mu$ through SSMLearn. We will then find the coefficients $b_{j,l}$ for the corresponding two-dimensional models on the SSMs, obtaining reduced-order, predictive models at some order of approximation $q$ in $\eta_1$ and $\mu$ \begin{equation} \begin{cases}\dot{\eta_1} &= \sum_{j+l\leq q}b_{j,l}\ \eta_1^j\ \mu^{\ l} \ + \  \mathcal{O}(q)\\ \dot{\mu}&=0\end{cases} \label{eq:parrom}\end{equation}

We wish to understand how the dynamics are modified by varying the input vector $\mathbf{u}$ to achieve the flexibility needed for a context-dependent decision-making task. We expect that varying context inputs $u_3 = c_1$ and $u_4= c_2$ would lead to significant changes in the domains of attraction of the two fixed points, whenever the sensory inputs have different signs.

We fix the sensory inputs to $s_1 = -0.15, \, s_2 = 0.15$ and parametrically vary the context inputs through the scalar $k$,  $c_1 = 1-k$, $c_2 = k$, $k \in [0,1]$. The number and stability of the fixed points remain unchanged for different values of $k$, but we are interested in how the the global phase space geometry changes with $k$, which depends on the locations of the boundaries of the domains of attraction of the two stable fixed points.

The 2D invariant manifold in the extended phase space $\mathbb{R}^{N+1}$ with coordinates $(\mathbf{y}, \bar{k})$, $\bar{k} = k -k_0$, $k_0 =0.5$, exists and is smooth in $\mathbf{y}$ and $k$, for every value of $k \in [0,1]$, due to the smooth dependence of the SSMs on the parameters, when the fixed points they are attached to are robust. Hence, we look for this invariant manifold as a graph over the extended spectral subspace $P$ through SSMLearn. We plot the slices of the extended phase space for different values of $k$ of the approximation of this manifold at order five in Fig. \ref{fig:kmflds} (upper left), and the whole invariant manifold in the extended phase space in Fig. \ref{fig:kmflds} (upper right). 
We find the corresponding reduced model at order nine as in Eq. \ref{eq:parrom}, with $\mu = \bar{k}$ and $q = 9$.

  We plot the parameter-dependent right-hand side in Fig. \ref{fig:kmflds} (upper center). We deduce from the 1D model the widths of the domains of attraction of the two stable fixed points in the reduced model by looking at the location of the unstable fixed point that acts as a $0$-dimensional boundary of such domains. These widths vary as expected for the network trained for the context-dependent task, since initial conditions around zero switch the fixed point they convergence to for $k=0.5$.

Next, we study the effects of varying sensory inputs on the system's dynamics.
First, we fix the context inputs at $c_1=0$, $c_2=1$ (context 2) and the select sensory input $s_2=0.036$, while we vary the sensory input $s_1$ in the range $[-0.036,0.036]$. Since the network must ignore the sensory clues coming from the input $s_1$, we expect that the network dynamics will look the same when $s_1$ varies in the considered range.

Since the fixed point configuration is unchanged in the considered range of $s_1$, we can look for the 2D invariant manifold (which exists and is smooth by the smooth dependence of unstable manifolds on parameters) in the extended phase space $\mathbb{R}^{N+1}$ with coordinates $(\mathbf{y}, \bar{s}_1)$, $\bar{s}_1 = s_1 -s_{1,0}$, $s_{1,0}=0$. The slices of this manifold for fixed $s_1$ correspond to the unstable manifolds of the unstable fixed points. We find an approximation of order five of the global manifold that yields a manifold fitting error of order $10^{-4}$.

In Fig. \ref{fig:kmflds} (bottom left), we plot the resulting 1D, parameter-dependent manifolds obtained for fixed values of the parameter $s_1$ in the above equation, while on the bottom left we plot the global invariant manifold in the extended phase space.
We also obtain the reduced-order model at order nine in $\eta_1$ and $s_1$ 
as in Eq. \ref{eq:parrom} with $\mu = s_1$ and $q = 9$
whose right-hand sides for the different values of $s_1$ are plotted in Fig.\ref{fig:kmflds} (bottom center).

As expected, the phase space geometry does not undergo significant changes when varying the sensory input $s_1$, if context two is selected. The parameter-dependent reduced-order model on the $s_1$-dependent unstable manifolds correctly captures this phenomenon. 

Finally, we study the consequences of parametric changes in sensory input 2 when context 2 is selected, and the first sensory input $s_1$ is kept fixed ($s_1=0.036$) in the range for $s_2$ selected by Mante et al. in \cite{mante_context-dependent_2013}, $s_2\in [-0.15,0.15]$.

\begin{figure}
\includegraphics[width=0.5\linewidth]{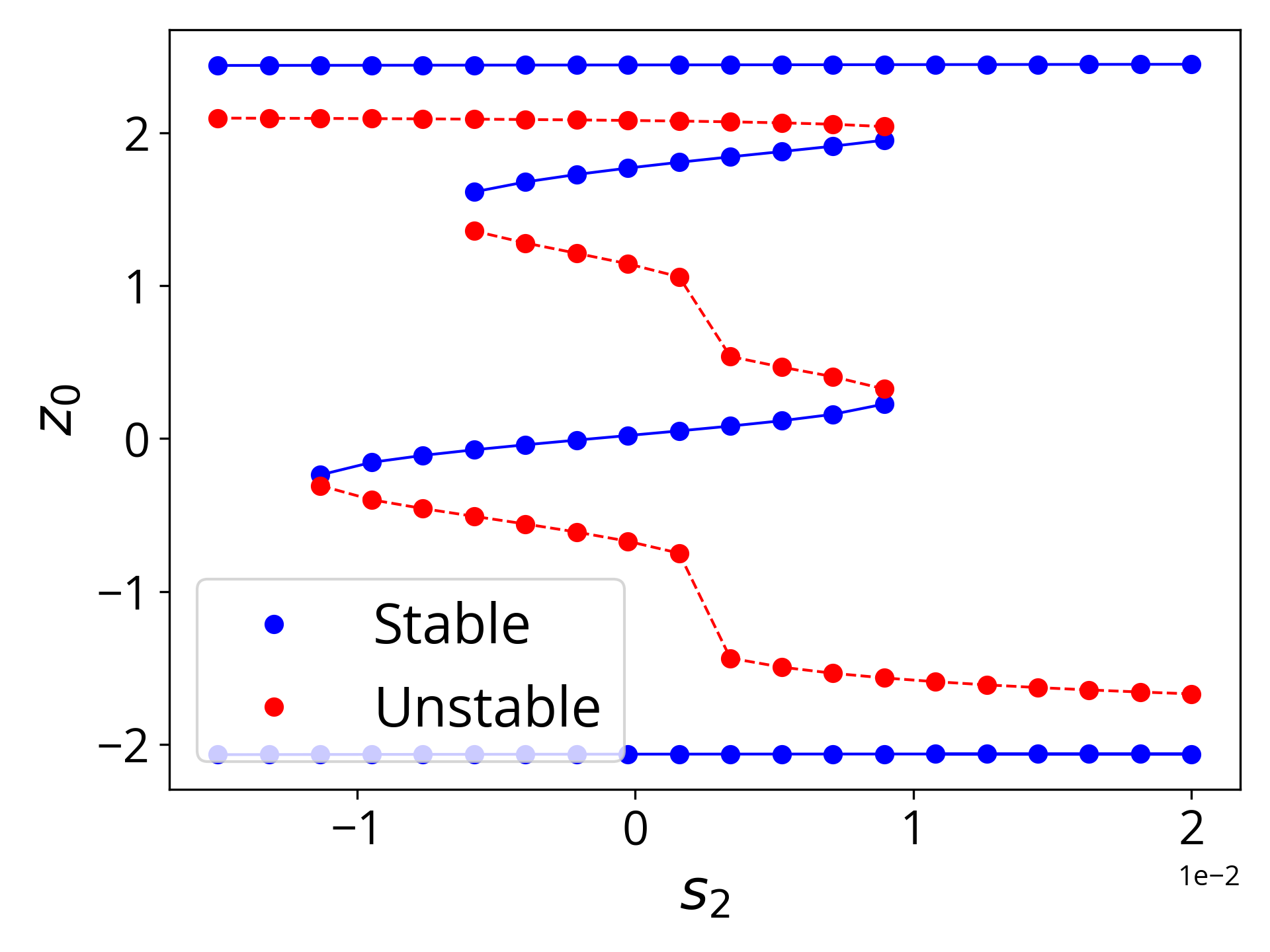}
    \includegraphics[width=0.45\linewidth]{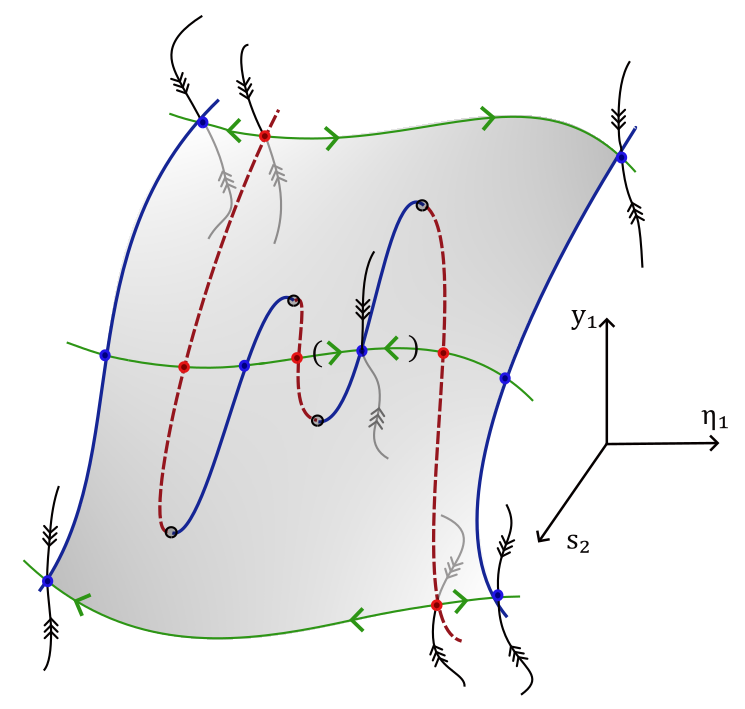}
    \centering
     
    \caption{(Left) Bifurcation diagram for the fixed points of the RNN, when sensory input 2 $s_2$ (selected by context input 2) is varied in the range $s_2 \in [-0.15, 0.15]$ and sensory input 1 is kept fixed at $s_1 = 0.036$. Four saddle-node bifurcations occur in the plotted critical parameter range $[-0.015, 0.020]$ around zero: two create a pair of fixed points, one stable and one unstable, and two annihilate one of such pairs. Through this cascade of saddle-node bifurcations, the system adapts by converging to a fixed point that yields the correct readout: negative when \( s_2 < 0 \), and positive when \( s_2 > 0 \). We plot a one-dimensional readout defined by $z_0 = \langle \mathbf{Y}, \mathbf{x}_0\rangle$. 
    (Right) Sketch of the robust, normally hyperbolic 2D slow manifold in the extended phase space, the global invariant manifold that carries the reduced dynamics of the RNN. When $s_2<0$ ($s_2>0)$, the negative-readout (positive-readout) fixed point has the largest domain of attraction, and trajectories converge to it along a branch of the unstable manifold of the unstable fixed point that acts as a boundary. When $s_2 = 0$, the stable fixed point around zero appeared via a saddle-node bifurcation attracts trajectories along its 1D slowest SSM (in round brackets).
 }
    \label{fig:bif_diag}
\end{figure}

\begin{figure}
   \centering
    \includegraphics[width=1\linewidth]{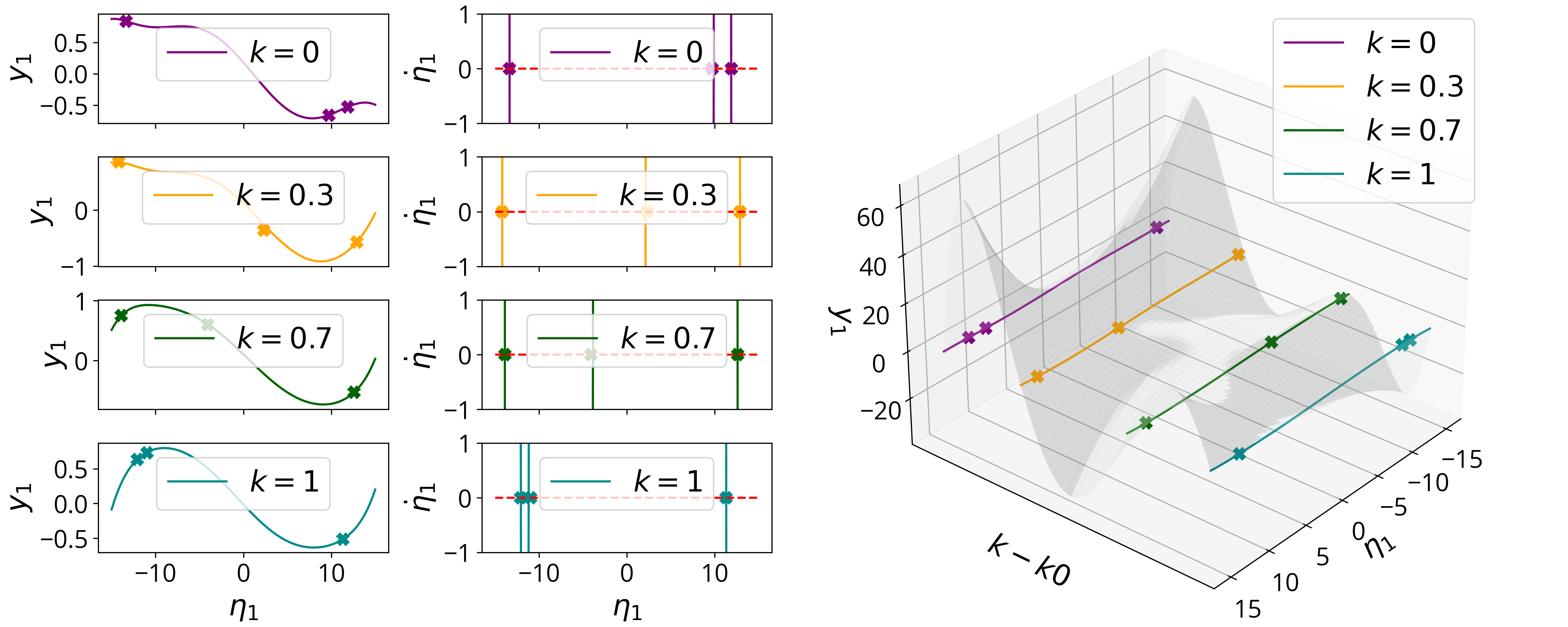}
    \centering
    \includegraphics[width=1\linewidth]{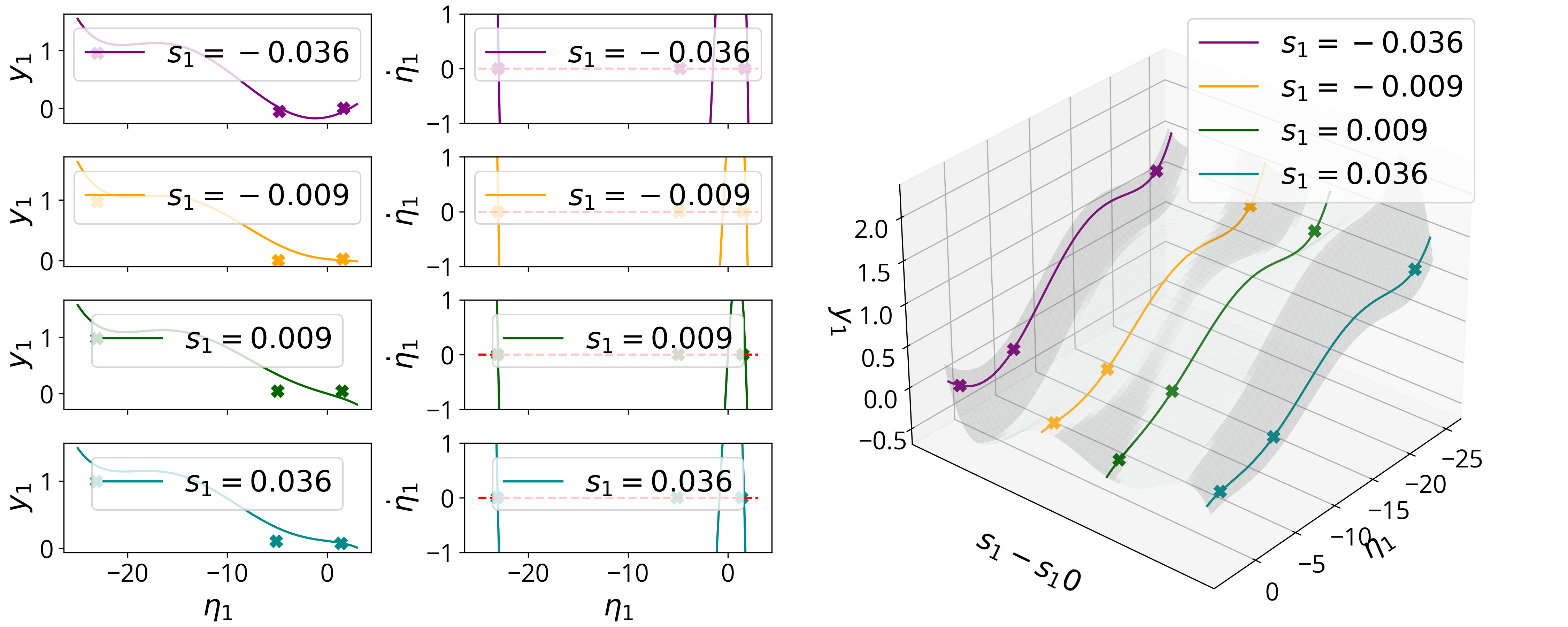}

    \caption{(Left) Parameter-dependent SSMs yielding the unstable manifolds for different values of $k$ (top), with $c_1= k $ and $c_2 = 1-k$, when $s_1 = -0.15$ and $s_2=0.15$, and for different values of $s_1$ (bottom), where $s_2 = 0.15$ and context 2 is selected. The manifolds are plotted in coordinates $(\eta_1, y_1,)$. The symbols $x$ indicate fixed points of the full model on the manifolds (two stable and one unstable in the middle, acting as boundary of the domains of attractions of the other two). (Center) Graphs of the right-hand sides of the parameter-dependent reduced-order models on the unstable manifolds for different values of $k$ (top) parametrizing the context inputs, $c_1 = k$ and $c_2 = 1-k$, and $s_1$ (sensory input 1). The symbols x indicate fixed points of the full-order model, to be compared to points where $\dot{\eta}_1 = 0$. We note how the domains of attraction of the two stable fixed points change parametrically with the context input $k$ as expected. When varying $s_1$, instead, the reduced dynamics do not vary, as expected (bottom). (Right) Visualization of the 2D slow manifolds in the extended phase space in coordinates $(\eta_1, \mu, y_1)$, where $\mu = \bar{k} = k-k_0$ (top) and $\mu = \bar{s}_1 = s_1-s_{1,0}$ (bottom).
    }
    \label{fig:kmflds}
\end{figure}

We find the system's fixed points for varying values of \( s_2 \) within the specified range using a numerical solver and observe that, even though the two stable fixed points corresponding to the two choices remain fixed, the number and location of the other fixed points change with the parameter. In particular, four saddle-node bifurcations occur within a small interval around $s_2 = 0$, as illustrated in the bifurcation diagram in Fig. ~\ref{fig:bif_diag} (right).

It is through this cascade of saddle-node bifurcations that the system adapts to favor the convergence to the stable fixed point that yields the correct readout (positive for positive values of $s_2$ and vice versa) along the branches of the 1D unstable manifold of an unstable fixed point.

If we take negative values of $s_2$, the unstable manifold that carries the reduced dynamics of the RNN is unique and normally hyperbolic. 
We show in Appendix~\ref{appA:paramSSM} that even when the associated unstable fixed point (the bottom-most in Fig.~\ref{fig:bif_diag}, right) disappears, a 1D invariant slow manifold persists throughout the entire parameter range and continues to carry the reduced dynamics of the network.

However, since the unstable fixed point to which we attach the SSM does not persist along the whole $s_2$ range, it is not possible to find a single parameter-dependent reduced model for the dynamics through SSMLearn. Nevertheless, for each value of $s_2$, we can find a 1D reduced model for the observed dynamics on the 1D slow SSM, and we can describe the network dynamics globally for the whole parameter range in the extended phase space (see Fig. \ref{fig:bif_diag}, left).

In particular, we already knew (see Fig. \ref{fig:unstablemfld}, bottom) that when $s_2$ is positive, the 1D unstable manifold of an unstable fixed point carries the decision-making reduced dynamics, leading convergence to the stable fixed point that corresponds to a positive readout. On the other hand, when $s_2$ is negative, the unstable manifold of the bottom-most unstable fixed point in Fig. \ref{fig:bif_diag} (right) arising from a saddle-node bifurcation when $s_2 \approx 0$ leads convergence to the fixed point corresponding to a negative value of the readout. These unstable manifolds correspond to line attractors when sensory inputs are on in \cite{mante_context-dependent_2013}.

When $s_2 \approx 0$, instead, the stable fixed point appeared in the critical $s_2$ interval around zero (where the four saddle-node bifurcations take place) attracts all the trajectories in a ball $B_{5\delta}$, $\delta = 0.1$ around zero, carrying a 1D slowest SSM (line attractor when sensory inputs are off in \cite{mante_context-dependent_2013}).

These results yield a comprehensive interpretation of the system’s behavior as governed by a global, 1D slow manifold structure that smoothly deforms for input parameter changes and whose internal dynamical structure adapts to consistently lead to convergence to the correct fixed point indicated by the inputs.



\section{Finite-Time Lyapunov Exponent (FTLE)}
\label{appC}
\begin{figure}
    \centering
    \includegraphics[width=1\linewidth]{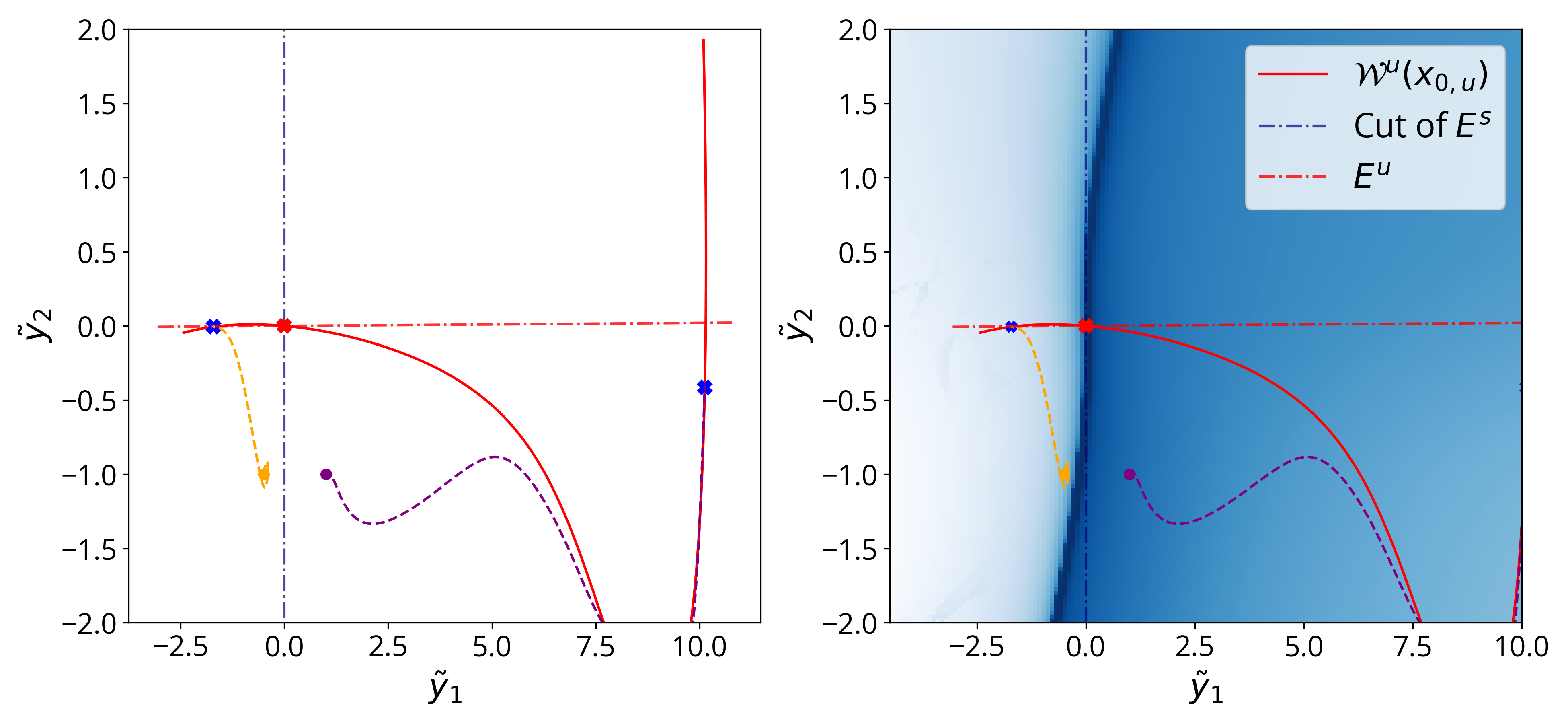}
    
    \caption{(Left) Projection of the dynamics around the unstable fixed point on the spectral subspace $E_2$. We plot the projected intersections of the unstable and stable subspaces ($E^S$ and $E^U$, respectively). (Right) $\mathrm{FTLE}$ calculated on the chosen two-dimensional plane $P=E_2$. Darker values in the picture on the left indicate higher values of the FTLE. Regions of high values of the FTLE indicate the estimated boundaries of the domains of attraction of the two fixed points. As expected, the cut of the stable manifold $\mathcal{W}^S$ on the plane $P$ is tangent to the stable subspace $E^S$. Trajectories with initial conditions belonging to the two different regions converge to different fixed points.}
    \label{fig:ftle2d}
\end{figure}
To estimate the sensitivity of RNN trajectories to initial conditions, we compute the \textit{Finite-Time Lyapunov Exponent (FTLE)}. Given the flow map \( F_{t_0}^t(x_0) \) for the system \( \dot{x} = f(x, t) \), the FTLE quantifies the maximal rate of trajectory separation:

\begin{equation}
\mathrm{FTLE}_{t_0}^t(x_0) = \frac{1}{2|t - t_0|} \log \lambda_{\max}(C_{t_0}^t(x_0)),
\end{equation}

where \( C_{t_0}^t = [DF_{t_0}^t(x_0)]^T DF_{t_0}^t(x_0) \) is the Cauchy-Green strain tensor (see \cite{TBarrier} for more details).

We evaluate the FTLE over a grid of initial conditions on the two-dimensional spectral subspace \( E_2 = \mathrm{span}\{\mathbf{e}_1, \mathbf{e}_2\} \), aligned via a linear change of basis around the unstable fixed point. The transformed system is:
    \begin{align} \dot{\mathbf{\tilde{y}}} = -\mathbf{\tilde{y}} - T \mathbf{x_0} + T \mathbf{W}\tanh(T^{-1}\mathbf{\tilde{y}}+\mathbf{x_0}) + T\mathbf{B\ u} \end{align} where \( \tilde{y} \in \mathbb{R}^N \) are the new coordinates and $T$ the transformation. Simulating the dynamics allows us to compute \( D_{\tilde{y}_1, \tilde{y}_2}F^t(\tilde{\mathbf{y}}_0) \) for any point $\tilde{\mathbf{y}}_0$ on the plane numerically by finite differences and, from it, the restricted FTLE.

FTLE ridges, local maxima of the FTLE field (see \cite{eberly}), indicate the boundaries between the domains of attraction of the two stable fixed points (see Fig.~\ref{fig:ftle2d}). These ridges correspond to intersections of the stable manifold \( W^S(x_0) \) of the unstable fixed point with the selected subspace.
\clearpage

\begin{center}
\textbf{\large Supplementary Methods}
\end{center}
\section{Persistence of the Slow Manifold for the Context-Dependent Decision-Making task}
\label{appA:paramSSM}
\begin{figure}
    \centering
    \includegraphics[width=\linewidth]{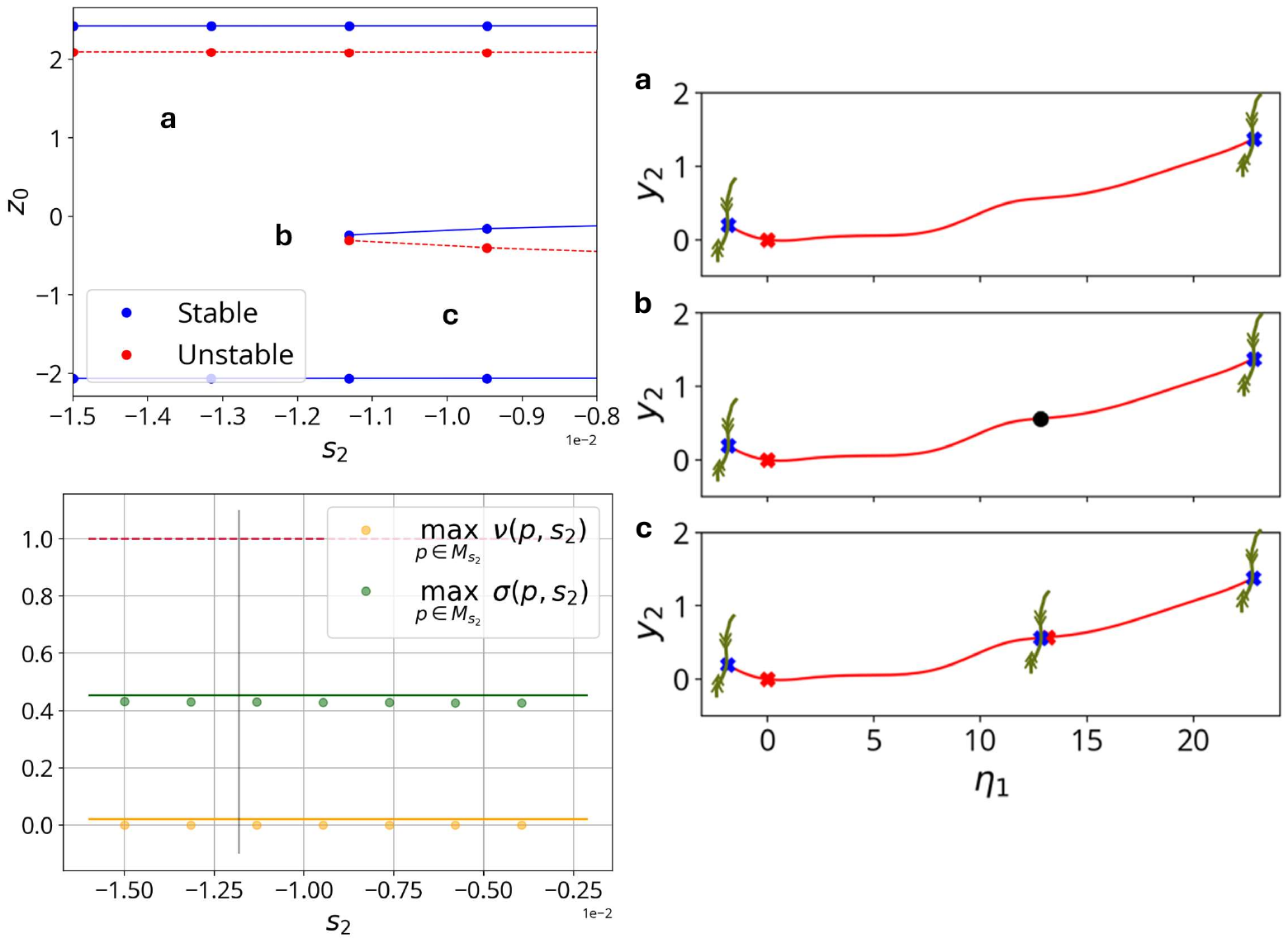}
    \caption{(Upper left and right) Bifurcation diagram in the restricted parameter range where the first saddle-node bifurcation happens along a 1D center manifold, and the corresponding normally attracting, unique invariant manifolds in cases (a), (b) and (c). This is the bifurcation responsible for the appearance of the unstable fixed point from which we construct the unstable manifold for positive values of $s_2$ (see Fig. \ref{fig:bif_diag} and Fig. \ref{fig:unstablemfld}, bottom).  (Bottom right) Supremum of type numbers for $\alpha$-limit sets of the invariant manifold over the $s_2$ range. Dots: numerical values; lines: upper bounds. Bounds below 1 indicate normal hyperbolicity throughout.}
    \label{fig:bif_diag_red} 
\end{figure}

In the context-dependent RNN of \cite{mante_context-dependent_2013}, the reduced dynamics lie on 1D unstable manifolds (when sensory inputs are on) or on the slowest SSM of a stable fixed point (when sensory inputs are off, $s_1 = s_2 = 0$) (see \ref{ss:cddmRNN}, \ref{fig:unstablemfld}). These manifolds coincide with the line attractor slow dynamics identified in \cite{mante_context-dependent_2013}.

To study the persistence of such a 1D slow manifold for the entire range of input parameters explored in the task, we consider $\epsilon$-parameter changes with $0<\epsilon  \ll 1$. Our goal is to determine whether the invariant manifold $\mathcal{W}_0$ of the system at $\epsilon = 0$ persists as a $C^{\rho}$-smooth manifold $\mathcal{W}_{\epsilon}$ that remains $\mathcal{O}(\epsilon)$-close to the original one.

Smooth dependence of SSMs on parameters generally suffices to state this persistence result, as long as the fixed points to which they are attached remain hyperbolic during parameter changes. In this case, though, when we vary the relevant sensory input for the task from negative to positive values ($s_1$ for context 1 and $s_2$ for context 2), we observe four saddle-node bifurcations shown in Fig. \ref{fig:bif_diag} (left). These bifurcations cause the unstable fixed points that we utilize to describe the reduced dynamics for both positive and negative sensory inputs through their unstable manifolds (see Fig. \ref{fig:bif_diag} and Fig. \ref{fig:unstablemfld}, bottom) to lose their hyperbolicity and disappear.

Nevertheless, a few simple observations let us establish a persistence result for this problem. We start by looking at negative relevant sensory input values ($s_2$ for context $c_2$) to focus on the parameter range where the first of such bifurcations happens, as shown in the bifurcation diagram in Fig. \ref{fig:bif_diag_red}.

In region (a), there are two stable and one unstable fixed points. To show there is a unique, attracting invariant manifold in the RNN phase space, we can use the asymptotic property of the unstable manifold of the latter, which is defined as the set of points that converge to the unstable fixed point asymptotically in backward time. Therefore, we can extend the local unstable manifold by forward-time integration of initial conditions along the unstable eigenspace in the proximity of the fixed point.

Fig. \ref{fig:bif_diag_red} (right, a) shows the results of forward integration, through which we can conclude that there exists a 1D attracting invariant manifold $\mathcal{W}_0$ composed of the two stable fixed points, the unstable fixed point, and the two branches of its unstable manifold.
This manifold is unique in class $C^0$ and, since the convergence to the stable fixed points along the unstable manifold must align with their slowest eigenvectors (the line spanned by $e_1$), $\mathcal{W}_0$ is also at least of class $C^1$. Moreover, it is compact with an invariant boundary, with the two stable fixed points serving as its closure. 

We can now invoke Fenichel’s results on the persistence of compact, overflowing-invariant manifolds (see \cite{Fenichel}, \cite{Wiggins}). These results require the manifold to be normally hyperbolic: the rate of contraction or expansion in directions transverse to the manifold must dominate the dynamics along the manifold itself by some integer $\rho$. This smoothing property of the flow ensures persistence and determines the degree of smoothness of the perturbed manifold $\mathcal{W}_\epsilon$.

While an exposition of the theory of Normally Hyperbolic Invariant Manifolds (NHIMs) is beyond our scope, we briefly state the main result: a compact, $C^r$, $\rho$-NHIM will persist under small $C^1$ $\mathcal{O}(\epsilon)$-perturbations as a nearby invariant manifold of differentiability class $C^{\mathrm{max\{r, \rho\}}}$ for the perturbed flow (see \cite{Fenichel, Wiggins} for rigorous statements and proofs).

We now wish to show that the invariant manifold $\mathcal{W}_0$ is an attracting NHIM. To do so, we need to verify the definition of normal hyperbolicity through the exponential rate of normal attraction and tangential contraction of the linearized dynamics along the manifold. The asymptotic properties of the linearized flow are quantified through the \emph{Lyapunov-type numbers}, for which we again refer to \cite{Fenichel} and \cite{Wiggins}. Since these numbers are defined by looking at asymptotic properties of trajectories on the manifold in backward time, to determine the normal hyperbolicity of a manifold, it is sufficient to calculate them on its $\alpha$-limit sets, i.e., the sets to which trajectories on the manifold converge in forward time. 

In this case, the only limit sets are the fixed points. In particular, we have that the invariant manifold $\mathcal{W}_0$ is $\rho$-normally attracting if, given the ordered eigenvalues of the linearized dynamics around a fixed point $p$, $\lambda_1 \geq \lambda_2 \geq \dots \geq \lambda_N$
    \begin{enumerate}
        \item The first Lyapunov type number quantifying the asymptotic rate of normal attraction  $\nu(p) = e ^ {\lambda_2} < 1$ for all the fixed points $p$ in $W_0$ ($\mathcal{W}_0$ attracts nearby trajectories)
        \item The second Lyapunov type number $\sigma(p) = \frac{\lambda_1}{\lambda_2}  < \frac 1 r$, for some integer $r>1$, for all the fixed points on $\mathcal{W}_0$ (the normal decay on $\mathcal{W}_0$ dominates over the tangential compression along $\mathcal{W}_0$ (if present))
    \end{enumerate}
As shown in Fig. \ref{fig:bif_diag_red} (bottom left), not only are the requirements on the Lyapunov type numbers satisfied, but there also exists a uniform bound along the whole considered parameter range. 

Therefore, by Fenichel's theorem, we can conclude that there exists a semi-open interval of parameters $\epsilon \in [ 0, \bar{\epsilon})$ (which parametrizes some semi-open interval for $s_2$ $[-0.15, s_2^*)$) along which this manifold exists, is at least $C^1$ and is normally attracting. We now wish to prove that we can extend this persistence result to include the parameter value $\bar{s}_2$ corresponding to the bifurcation point (b) in Fig. \ref{fig:bif_diag_red}, which is the only value where we cannot use the uniqueness and asymptotic properties of the unstable manifold to state the existence, uniqueness and smoothness result for the invariant manifold.

The uniform boundedness of the Lyapunov type numbers below one is a necessary condition to conclude this result (see \cite{kopell}). As stated in \cite{chiconeNHIMS}, though, it is not sufficient to determine the existence of a $C^1$-smooth manifold at the boundary of the parameter range. In fact, this manifold might lose normal hyperbolicity by losing $C^1$ smoothness. 

The complete argument goes as follows.
Let $A = \{\epsilon \in [0,\bar{\epsilon}]: \forall \epsilon' \in [0,\epsilon]$, there exists a NHIM$\mathcal{W}_{\epsilon'}$$\}$. If $A \neq \emptyset $, $A$ open and $A$ closed, then $\bar{\epsilon} \in A$, and the result follows. We know that $A$ is nonempty, since $\mathcal{W}(0)$ is a NHIM, and the fact that $A$ is open follows from Fenichel's theorem. To show that $A$ is closed, we need to show that if $\epsilon^*$ is the supremum of $A$, then $\epsilon^*\in A$. In other words, we need to find an argument for the existence of a $C^1$ invariant manifold for $\epsilon^*$. 

In our setting, the only value in $A$ we need to check is the bifurcation value $\epsilon^* = \bar{\epsilon}$. Since we know from calculating the Jacobian that the saddle-node bifurcation that gives rise to the pair of unstable and stable fixed points has a one-dimensional center manifold, and the center manifold is $C^r$-smooth for $r \in \mathbb{N}$, the existence of a $C^1$ invariant manifold is guaranteed at the bifurcation value. The manifold is normally hyperbolic, as shown in Fig. \ref{fig:bif_diag_red}.

To show the existence of an attracting slow manifold in the extended phase space, describing the reduced dynamics of the context-dependent decision-making RNN for the whole considered $s_2$ range as in Fig. \ref{fig:bif_diag} (right), it remains to show that these normally hyperbolic manifolds we proved to exist are also unique. 

We already know that in region (a), the manifold is unique by the uniqueness of the unstable manifold. In region (c), we can similarly invoke the definition of unstable manifold and continue the 1D branches emanating from each unstable fixed point by forward integration of trajectories along their unstable direction. The resulting invariant manifold, composed of three stable fixed points, two unstable fixed points, and the four branches of their unstable manifolds, is therefore also unique. The uniqueness at the bifurcation point follows from the uniqueness of the limit.

\section{Normal Hyperbolicity of the Heteroclinic Structure for the Memory-Pro Task}
\label{appA:wm}
We aim to establish the normal hyperbolicity of the heteroclinic manifold \( \mathcal{W}_h \) in the 2D reduced-order model of the multitasking RNN presented in \cite{ringattractor}, which performs the Memory-Pro task.

The manifold \(\mathcal{W}_h \) consists of a stable fixed point, an unstable fixed point, and the two branches of the unstable manifold connecting them. This set forms a compact, boundaryless, and invariant manifold by construction.

As shown in Figure~\ref{fig:slowring}, trajectories initially converge onto the two-dimensional unstable manifold that contains \( \mathcal{W}_h \), and subsequently contract rapidly onto \( \mathcal{W}_h \) itself, ultimately reaching the stable fixed point. To verify normal hyperbolicity, we evaluate the Fenichel-type numbers \( \nu(p) \) and \( \sigma(p) \) for all \( p \in \mathcal{W}_h \), characterizing the asymptotic behavior as \( t \to +\infty \).

The only limit sets within \( \mathcal{W}_h \) are the two fixed points. The unstable one has a 1D unstable manifold and N-1 decaying directions in the linearized dynamics. The linearized dynamics around the stable fixed point have the first two eigenvalues satisfying \( \mathrm{Re}\lambda_{2,s} < \mathrm{Re}\lambda_{1,s} < 0 \). The eigenvector corresponding to \( \lambda_2 \) is tangent to \( \mathcal{W}_h \), while the one associated with \( \lambda_1 \) is transverse. We can calculate:
\begin{itemize}
    \item $\max_{p \in \mathcal{W}_h}\nu(p) = \max {\{e^{\lambda_2,s},\ e^{\lambda_2,u}\}} = 0.0156 < 1$ (normal attraction rate)
    \item $\max_{p \in \mathcal{W}_h}\sigma(p) = \frac{\lambda_{1,s}}{\lambda_{2,s}} = \frac{-0.0496}{-5.0213} < \frac{1}{101}$
\end{itemize}
Hence, $\mathcal{W}_h$ is $\rho$-normally attracting with $\rho = 101$.

\section{Nonautonomous Spectral Submanifolds (SSMs)}
\label{appB}

SSM theory extends to nonautonomous systems with periodic, quasiperiodic, weak aperiodic, slow aperiodic and even random forcing. We briefly summarize the key results and refer the reader to \cite{firstSSM, aperiodicSSM, randomSSM} for formal statements, assumptions, and proofs.

For periodically or quasiperiodically forced systems of the form
\begin{equation}
    \dot{x} = Ax + f_0(x) + \epsilon f_1(x, \Omega t; \epsilon), \quad x \in \mathbb{R}^N, \ \Omega \in \mathbb{R}^k,
\end{equation}
with smooth nonlinearities $f_0 = \mathcal{O}(\|x\|^2)$ and $f_1$ periodic in each $\Omega_i t$, fixed points of the unforced system perturb into periodic (or quasiperiodic) solutions known as \textit{Nonlinear Normal Modes (NNMs)}. Invariant manifolds dependent on time, non-autonomous SSMs, exist in their vicinity under non-resonance conditions \cite{firstSSM}.

In the RNN analysis, the only time-dependent input is bounded aperiodic forcing, in the form of drawings from a bounded Gaussian distribution. For this setting, we consider systems of the form
\begin{equation}
    \dot{x} = Ax + f_0(x) + f_1(x, t),
\end{equation}
where $f_1$ is uniformly bounded in time but not necessarily continuous. If the unforced system has a hyperbolic fixed point, this point perturbs into a uniformly bounded, compact trajectory called an \textit{anchor trajectory} under specific uniform boundedness conditions on the time-dependent forcing $f_1$ \cite{aperiodicSSM}. SSMs can then be constructed along this trajectory, provided a spectral subspace $E$ is $\rho$-normally hyperbolic and satisfies nonresonance conditions \begin{equation}\sum_{j=1}^Nm_j \mathrm{Re}\lambda_j \neq \mathrm{Re}\lambda_k \ \ \forall k \in \{1, \dots, n\}, \quad m_j \in \mathbb{N}, \ \ \sum_{i=1}^N m_j \geq2.\end{equation}

The key result is that under these assumptions, a unique time-dependent SSM $\mathcal{W}_E(x^*(t))$ exists, is tangent to $E$ at each time step, and can be computed as the uniformly bounded solution of a time-dependent nonlinear PDE, even for larger forcing amplitudes than predicted by theory. These manifolds enable model reduction and capture the slow nonlinear dynamics of recurrent neural networks (RNNs) under single realizations of bounded Gaussian noise in the form of weak aperiodic forcing.

Formal expansions for both anchor trajectories and time-dependent SSMs are derived in \cite{aperiodicSSM}. In our work, specific expansions for RNN dynamics are presented in ~\ref{appB:anchor} and~\ref{appB:SSMtdep}.

\subsection{Anchor Trajectory Calculations}
\label{appB:anchor}
First, we want to compute anchor trajectories $\mathbf{x}^*(t)$ perturbing from fixed points when small-amplitude discontinuous forcing is added to the equations of the autonomous RNNs.

We have the rescaled equations of motions for the Vanilla RNNs \begin{align}
    \mathbf{x}'= \frac{d\mathbf{x}}{dt'} = \frac{1}{\tau} \frac{d\mathbf{x}}{dt} = \frac{1}{\tau} \Big(-\mathbf{x} + \mathbf{W}\,\mathbf{r(x)} +  \mathbf{B}\,\mathbf{u}\Big)
\end{align}

We change coordinates to $\mathbf{y}= \mathbf{x}-\mathbf{x_0}$ so that the fixed points we want to attach SSMs to are in zero. In the following, we change notation from $\mathbf{y}'$ to $\mathbf{\dot{y}}$ to indicate the derivative with respect to the rescaled time $\frac{t}{\tau}$, and we denote such rescaled time with $t$. Hence,
    \begin{align}\begin{split}\label{eq:genAp} \dot {\mathbf{y}} &= \mathbf{A}\,(\mathbf{y}-\mathbf{x_0}) +\mathbf{f(y},t),\\          \mathbf{f(y},t)&=\mathbf{f}_0({\mathbf{y}})+\epsilon \mathbf{f}_1(\mathbf{y},t) + \mathbf{Bu},   \quad  \mathbf{f}_0(\mathbf{y})=O(|\mathbf{y}|^2), \end{split}\end{align}
with  
    \begin{align*} \mathbf{A} &=  -\mathbf{I} + \begin{pmatrix} \frac{w_{11}}{\cosh^2(x_{1,0})} & \dots &\dots & \frac{w_{1N}}{\cosh^2(x_{N,0})} \\ \vdots &\vdots&\vdots &\vdots \\ \frac{w_{N1}}{\cosh^2(x_{1,0})} &\dots&  \dots &\frac{w_{NN}}{\cosh^2(x_{N,0})} &\end{pmatrix}, \\ \mathbf{f}_0(\mathbf{y}) &= \mathbf{W\tanh}(\mathbf{y}-\mathbf{x_0}) -  \begin{pmatrix} \frac{w_{11}}{\cosh^2(x_{1,0})} & \dots &\dots & \frac{w_{1N}}{\cosh^2(x_{N,0})} \\ \vdots &\vdots&\vdots &\vdots \\ \frac{w_{N1}}{\cosh^2(x_{1,0})} &\dots&  \dots &\frac{w_{NN}}{\cosh^2(x_{N,0})} &\end{pmatrix}(\mathbf{y}-\mathbf{x}_0),\\\mathbf{f}_1(\mathbf{y},t) &= \mathbf{f}_1(t) =  \mathbf{\sigma}(t),\end{align*} 
    and $\mathbf{\sigma}(t)$ is drawn at each time step from a bounded standard normal distribution. 
    
    
 


We use the formal expansions in \cite{aperiodicSSM} to approximate anchor trajectories perturbing from fixed points through $
\mathbf{y}^*(t) =\sum_{\nu=1}^M \epsilon^{\nu}\mathbf{y}_{\nu}(t)+o(|\epsilon\mathbf{f}_1|^M_U)
$, assuming it exists. 

We calculate the expansion terms $\mathbf{y}_{\nu}(t)$ up to order $M=3$, and perform (discrete) numerical integration to obtain these contributions numerically at each time step. 
For the first-order term, given $\mathbf{A}$ with $d_u$ unstable eigenvalues, we have 

\begin{align}\begin{split} \mathbf{y}_1(t) &= \int_{-\infty}^t\mathbf{T}e^{\begin{pmatrix}
    \mathbf{0} & \mathbf{0} \\ \mathbf{0} & \mathbf{A}_s \end{pmatrix}(t-\tau)}\mathbf{T}^{-1}\mathbf{f}_1(\tau)\,d\tau \ -\int^{+\infty}_t\mathbf{T}e^{\begin{pmatrix}
    \mathbf{A}_u & \mathbf{0} \\ \mathbf{0} & \mathbf{0} \end{pmatrix}(t-\tau)}\mathbf{T}^{-1}\mathbf{f}_1(\tau)\,d\tau , 
\\ \mathbf{A} &= \mathbf{T}\begin{pmatrix}
    \mathbf{A}_u & \mathbf{0} \\ \mathbf{0} & \mathbf{A}_s 
\end{pmatrix}\mathbf{T}^{-1} \quad\mathbf{A}_u = \begin{pmatrix}
    \lambda_1  & \mathbf{0} & \dots \\ \mathbf{0} & \ddots & \vdots \\ \vdots & \dots & \lambda_{d_u} \end{pmatrix} \quad\mathbf{A}_u = \begin{pmatrix}
    \lambda_{d_u+1}  & \mathbf{0} & \dots \\ \mathbf{0} & \ddots & \vdots \\ \vdots & \dots & \lambda_{N} 
\end{pmatrix},
\label{eq:anchoro1}\end{split}\end{align}
We write the expressions for the second and third-order contributions
        \begin{align*}\mathbf{y}_2(t) &= \int_{-\infty}^t\mathbf{T}e^{\begin{pmatrix}
    \mathbf{0} & \mathbf{0} \\ \mathbf{0} & \mathbf{A}_s \end{pmatrix}(t-\tau)}\mathbf{T}^{-1}[\,\frac{1}{2}\partial^2_{\mathbf{y}}\mathbf{f_0}(\mathbf{0})\otimes \mathbf{y}_1(\tau)\otimes \mathbf{y}_1(\tau)] \, d\tau \\ &-\int^{+\infty}_t\mathbf{T}e^{\begin{pmatrix}
    \mathbf{A}_u & \mathbf{0} \\ \mathbf{0} & \mathbf{0} \end{pmatrix}(t-\tau)}\mathbf{T}^{-1}[\,\frac{1}{2}\partial^2_{\mathbf{y}}\mathbf{f}_0(\mathbf{0})\otimes \mathbf{y}_1(\tau)\otimes \mathbf{y}_1(\tau)] \, d\tau \\ 
    \mathbf{y}_3(t) &= \int^t_{-\infty} \mathbf{T}e^{\begin{pmatrix}
    \mathbf{0} & \mathbf{0} \\ \mathbf{0} & \mathbf{A}_s \end{pmatrix}(t-\tau)}\mathbf{T}^{-1}\Big[\frac{1}{6} \partial_{\mathbf{y}}^3\mathbf{f}_0(\mathbf{0})\otimes \mathbf{y}_1(\tau)\otimes \mathbf{y}_1(\tau)\otimes \mathbf{y}_1(\tau) \\ &+ \frac{1}{2}D^2_{\mathbf{y}} \mathbf{f}_0(\mathbf{0})\otimes \mathbf{y}_1(\tau) \otimes \mathbf{y}_2(\tau)\Big]d\tau \\ &-  \int_t^{+\infty} \mathbf{T}e^{\begin{pmatrix}
    \mathbf{A}_u & \mathbf{0} \\ \mathbf{0} & \mathbf{0} \end{pmatrix}(t-\tau)}\mathbf{T}^{-1}\Big[\frac{1}{6} \partial_{\mathbf{y}}^3\mathbf{f}_0(\mathbf{0})\otimes \mathbf{y}_1(\tau)\otimes \mathbf{y}_1(\tau)\otimes \mathbf{y}_1(\tau) \\&+ \frac{1}{2}D^2_{\mathbf{y}} \mathbf{f}_0(\mathbf{0})\otimes \mathbf{y}_1(\tau) \otimes \mathbf{y}_2(\tau)\Big]d\tau,
\end{align*} where 
    \begin{align*}\mathbf{f}_0(\mathbf{y}) &= \begin{pmatrix}w_{11}\tanh(y_1-x_{1,0})+ ...  - (\frac{w_{11}}{\cosh^2(x_{1,0})}(y_1-x_{1,0})+ ... +\frac{w_{1N}}{\cosh^2(x_{N,0})}(y_N-x_{N,0})) \\ \vdots \\w_{N1}\tanh(y_1-x_{1,0})+ ... - (\frac{w_{N1}}{\cosh^2(x_{1,0})}(y_1-x_{1,0})+ ... +\frac{w_{NN}}{\cosh^2(x_{N,0}}(y_N-x_{N,0}))\end{pmatrix} ,\\ \partial_\mathbf{y}\mathbf{f}_0(\mathbf{x}_0) &= \mathbf{0},\\     
\partial_\mathbf{y}\mathbf{f}_0(\mathbf{y}) &= \begin{pmatrix} \frac{w_{11}}{\cosh^2(y_1-x_{1,0})} & \dots &\frac{w_{1N}}{\cosh^2(y_N-x_{N,0})} \\ \vdots & \vdots&\vdots \\\frac{w_{11}}{\cosh^2(y_1-x_{1,0})}&\dots &\frac{w_{1N}}{\cosh^2(y_N-x_{N,0})}\end{pmatrix}-\begin{pmatrix} \frac{w_{11}}{\cosh^2(x^0_1)} & \dots &\frac{w_{1N}}{\cosh^2(x^0_N)} \\ \vdots & \vdots&\vdots \\\frac{w_{11}}{\cosh^2(x^0_1)}&\dots &\frac{w_{1N}}{\cosh^2(x^0_N)}\end{pmatrix}.\end{align*}
       Note that $\partial^2_\mathbf{y}\mathbf{f}_0(\mathbf{y})$ is a 3-tensor wherein each vector entry is a Hessian matrix:
        \begin{align*}\partial^2_\mathbf{y}\mathbf{f}_0(\mathbf{y}) = \begin{pmatrix} \begin{bmatrix}-\frac{2w_{11}\sinh(y_1-x_{1,0})}{\cosh^3(y_1-x_{1,0})} & & \\ & \ddots & \\ & & -\frac{2w_{1N}\sinh(y_N-x_{N,0})}{\cosh^3(y_N-x_{N,0})}\end{bmatrix} \\\vdots \\\begin{bmatrix}-\frac{2w_{N1}\sinh(y_1-x_{1,0})}{\cosh^3(y_1-x_{1,0})} & & \\ & \ddots & \\ & & -\frac{2w_{NN}\sinh(y_N-x_{N,0})}{\cosh^3(y_N-x_{N,0})}\end{bmatrix} \end{pmatrix} \end{align*}
   Finally, $\partial^3_\mathbf{y}\mathbf{f}_0(\mathbf{y})$ is a 4-tensor (each $N$-vector entry is a $N\times N$ matrix of row $N$-vectors):

    \begin{align*}\partial^3_\mathbf{y}\mathbf{f}_0(\mathbf{y}) = \begin{pmatrix} \begin{bmatrix} \begin{bmatrix} w_{11}f^3(y_1-x_{1,0}) & 0 &\dots &0\end{bmatrix} & & \\ & \ddots & \\ & & \begin{bmatrix} 0 & \dots &0 &w_{1N}f^3(y_N-x_{N,0})\end{bmatrix}\end{bmatrix} \\\vdots \\\begin{bmatrix} \begin{bmatrix} w_{N1}f^3(y_1-x_{1,0}) & 0 &\dots &0\end{bmatrix} & & \\ & \ddots & \\ & & \begin{bmatrix}0 & \dots & 0 &w_{NN}f^3(y_N-x_{N,0}) \end{bmatrix}\end{bmatrix} \end{pmatrix} \end{align*}

    where
    
    $$
    f^3(y_i-x_{i,0}) = \frac{4\sinh^2(y_i-x_{i,0})}{\cosh^4(y_i-x_{i,0})} - \frac{2}{\cosh^4(y_i-x_{i,0})} = \frac{4\sinh^2(y_i-x_{i,0})-2}{\cosh^4(y_i-x_{i,0})}
    $$
    
\subsection{Nonautonomous SSM Calculations}
   \label{appB:SSMtdep}
We calculate the aperiodic, time-dependent SSMs, $\mathcal{W}(E,t)$, attached to the anchor trajectories obtained above $\mathbf{y}^*(t) = \mathbf{x}^*(t)-\mathbf{x}_0$, perturbing from the autonomous SSMs, $\mathcal{W}(E)$, tangent to a spectral subspace $E=E_1\oplus \dots \oplus E_k$. 

Let $\mathbf{P} = [\mathbf{e}_1\ \dots \ \mathbf{e}_N] \in \mathbb{C}^N$ be the matrix containing the complex eigenvectors corresponding to the ordered eigenvalues of $A$, and change coordinates to \begin{align*}\begin{pmatrix} \mathbf{u} \\ \mathbf{v}\end{pmatrix}= \mathbf{P}^{-1}(\mathbf{x}-\mathbf{x}^*_{\epsilon})\quad \implies \quad \begin{pmatrix} \dot{\mathbf{u}} \\ \dot{\mathbf{v}}\end{pmatrix}= \begin{pmatrix} \mathbf{\Lambda} & 0 \\ 0 &{\mathbf{A_v}}\end{pmatrix}\begin{pmatrix} u \\ \mathbf{v}\end{pmatrix}+ \hat{f}(u, \mathbf{v},\epsilon; t)\ \end{align*}
with\begin{align*}
    \mathbf{u} &\in \mathbb{R}^d, \ \mathbf{v} \in \mathbb{R}^{N-d}, \ \quad \mathbf{\Lambda} = \begin{pmatrix}
    \lambda_1 &&\\ & \ddots & \\
     && \lambda_d
\end{pmatrix}, \quad \mathbf{A_v} = \begin{pmatrix}
    \lambda_{d+1} &&\\ & \ddots & \\
     && \lambda_N
\end{pmatrix}, \\  \hat{f}&(\mathbf{u}, \mathbf{v},\epsilon; t) = \mathbf{P}^{-1}\Big[\mathbf{f_0}\Big(\mathbf{x}^*_{\epsilon}+\mathbf{P}\begin{pmatrix} \mathbf{u} \\ \mathbf{v}\end{pmatrix}\Big) + \mathbf{A}\mathbf{x}^*_{\epsilon}+\mathbf{A}\mathbf{\dot{x}}^*_{\epsilon} + \epsilon\, \mathbf{f_1}(t)\Big] 
\end{align*} 

In the case of the decision-making RNN, the time-dependent SSM, $\mathcal{W}(E_1,t)$, perturbs from the autonomous slowest SSM $\mathcal{W}(E_1)$ attached to the stable fixed point in the pre-stimulus period. 

The SSM, $\mathcal{W}(E_1,t)$, admits a formal asymptotic expansion    
        $$
        \mathcal{W}(E_1,t) = \Big\{(u,\mathbf{v}) \in U \subset\mathbb{R}\times\mathbb{R}^{N-1} \,;\, \mathbf{v} = \mathbf{h}_{\epsilon}(u,t) = \sum_{|(k,p)|\geq 1} ^M\mathbf{h}^{kp}(t)u^k\epsilon^p + o(|u|^q, \epsilon^{M-q})\Big\}.
        $$
According to \cite{aperiodicSSM}, we have \begin{itemize}
    \item $\mathbf{h}^{0,p} = 0$.
    \item  $\mathbf{h}^{k,0}$  are the coefficients of the unperturbed manifold $\mathcal{W}(E_1)$ expansion over the space spanned by ${\mathbf{e}_1}$.
    \item The other coefficients  $\mathbf{h}^{kp}$ can be obtained through the integral formula (\cite{aperiodicSSM}) $$
        \mathbf{h}^{kp}= \int_{0}^{\infty} \mathbf{G}_k(t-s) \,\mathbf{M}^{kp}(s, \mathbf{h}^{jm}(s)) \,ds
        $$
        and $\mathbf{G}_k(t) = e^{(\mathbf{A_v}-k\,\lambda_1\mathbf{I}_{N-1})t}$.
\end{itemize}
We stop at order $M = 2$, since first-order contributions already yield a good approximation of the time-dependent manifold. We have that \begin{align*}\mathcal{W}(E_1,t) = \Big\{(u,\mathbf{v}) \in U \subset \mathbb{R}^N \,;\, \mathbf{v} = \mathbf{h}_{\epsilon}(u,t) &= \sum_{|(k,p)|\geq 1} ^2\mathbf{h}^{kp}(t)u^k\epsilon^p + o(|u|^q, \epsilon^{2-q})\\ &= \mathbf{h}^{20}u^2 + \mathbf{h}^{11}u \ \epsilon + o( |u|^q, \epsilon^{2-q})\Big\},\end{align*} and
     \begin{itemize}
         \item The autonomous coefficients:\begin{itemize}
             \item  $\mathbf{h}^{10} = \mathbf{0}$.
             \item $\mathbf{h}^{20} = - \mathbf{A}_2^{-1}M^{20}(\mathbf{0}), \quad M^{20}(\mathbf{0}) =\Big[P^{-1}\frac{\partial^2\mathbf{f}_0}{\partial u^2}(\mathbf{x}_0)P\begin{pmatrix}1\\0\end{pmatrix}P\begin{pmatrix}1\\0\end{pmatrix}\Big]\Big|_{\mathbf{v}}$.
         \end{itemize}
         \item $k = 1, p = 1$ means the anchor trajectory needs to be inserted at first order, and $j=1$, $m=0$: 
    
    $$
    M^{11}(t,\mathbf{0})= \Big[P^{-1}\frac{\partial^2}{\partial^2 \mathbf{y}}\mathbf{f_0(x_0)}\otimes {\mathbf{y}}_1(t)\otimes P \begin{pmatrix}1 \\ 0\end{pmatrix}\Big]\Big|_{\mathbf{v}}.
    $$
     \end{itemize} 



\end{document}